\documentclass{amsart}

\usepackage{amssymb,amsmath,amsthm}
\usepackage[english]{babel}
\usepackage{enumerate}

\newtheorem{theorem}{Theorem}[section]
\newtheorem{proposition}[theorem]{Proposition}
\newtheorem{lemma}[theorem]{Lemma}
\newtheorem{corollary}[theorem]{Corollary}
\newtheorem{thmx}{Theorem}

\theoremstyle{definition}
\newtheorem{definition}[theorem]{Definition}
\newtheorem{example}[theorem]{Example}
\newtheorem{setting}[theorem]{Setting}
\theoremstyle{remark}
\newtheorem{remark}[theorem]{Remark}

\numberwithin{equation}{section}
\newcommand{\Z}{\mathbb{Z}}
\newcommand{\Q}{\mathbb{Q}}
\newcommand{\C}{\mathbb{C}}
\newcommand{\R}{\mathbb{R}}
\newcommand{\PP}{\mathbb{P}}
\newcommand{\AAA}{\mathbb{A}}
\newcommand{\F}{\mathbb{F}}
\newcommand{\cX}{\mathcal{X}}
\newcommand{\cU}{\mathcal{U}}
\newcommand{\cV}{\mathcal{V}}
\newcommand{\cO}{\mathcal{O}}
\newcommand{\cB}{\mathcal{B}}
\newcommand{\id}{\mathrm{id}}
\newcommand{\pr}{\mathrm{pr}}
\DeclareMathOperator{\Alb}{Alb}
\DeclareMathOperator{\End}{End}
\DeclareMathOperator{\Hom}{Hom}
\DeclareMathOperator{\Aut}{Aut}
\DeclareMathOperator{\Prym}{Prym}
\DeclareMathOperator{\Mor}{Mor}
\DeclareMathOperator{\Gal}{Gal}
\DeclareMathOperator{\rk}{rk}
\DeclareMathOperator{\Spec}{Spec}

\begin{document}

\title[Equivariant homomorphisms and twists of Albanese varieties]{Equivariant homomorphisms and Mordell--Weil groups of twists of Albanese varieties of cyclic multiple planes}

\author{Sajad Salami}
\address{Sajad Salami\\
Institute of Mathematics and Statistics,\\
 State University of Rio de Janeiro, Rio de Janeiro, Brazil}
\email{sajad.salami@ime.uerj.br}
\urladdr{https://sites.google.com/a/ime.uerj.br/sajadsalami/home}


\begin{abstract}
Let $X'\to X$ be a Galois cover of smooth projective varieties over a field $k$ with group $G$, and let $G$ act on an abelian variety $A$ over $k$. We prove an exact sequence relating the Mordell--Weil group of the twist of $A$ by $k(X')/k(X)$ to the group of $G$-equivariant homomorphisms from the Albanese variety of $X'$ to $A$; the remaining obstruction is killed by $|G|$. This corrects structure theorems in two earlier papers of the author and in a corrigendum to the first, where the equivariance condition was omitted. For the twists of Albanese varieties of cyclic multiple planes considered in those papers, we determine the rank and the torsion subgroup. Over number fields we describe the rank through the Galois action on endomorphisms and prove a specialization result. When $F(x,y)=0$ is the dual of the Fermat cubic, the Mordell curves $y^2=x^3+16F(a,b)$ have rank at least $3$ over $\Q$ for all $(a,b)\in\Q^2$ outside a thin set.
\end{abstract}

\maketitle

\section{Introduction}\label{sec:intro}

\subsection{Background}
Twisting a constant abelian variety by a finite extension of a function field is a classical way of producing abelian varieties with large Mordell--Weil rank; see for instance Tate and Shafarevich~\cite{TS} and Ulmer~\cite{Ul}. In~\cite{S19} the author considered the following situation, with the aim of extending a result of Hazama~\cite{Ha} on quadratic twists to cyclic covers of arbitrary degree. Let $k$ be a global field of characteristic prime to $n$ containing a primitive $n$-th root of unity, let $\pi\colon X'\to X$ be a cyclic cover of degree $n$ of varieties over $k$ with $X'(k)\neq\emptyset$, let $L=k(X')$ and $K=k(X)$, and let $\tilde A$ be the twist of an abelian variety $A$ over $k$ by $L/K$ and an automorphism $\sigma$ of $A$ of order $n$. Then \cite[Theorem~1.1]{S19} asserts that
\[
\tilde A(K)\cong \Hom_k(\Prym_{X'/X},A)\oplus A[n](k),
\]
and that $\rk\tilde A(K)\ge m\cdot\rk\End_k(A)$ whenever $\Prym_{X'/X}$ is $k$-isogenous to $A^m\times B$, with conditions on $B$. In the corrigendum~\cite{S21} the varieties were restricted to be smooth and projective, the conditions on the dimension of $B$ were removed, and the inequality was replaced by the equality $\rk\tilde A(K)=m\cdot\rk\End_k(A)$; this amended statement is also \cite[Theorem~3.3]{S20}. In the same way, \cite{S21} replaced the inequality in \cite[Theorem~1.2]{S19}, which applies the structure theorem to the superelliptic curves $y^n=f(x)$, by an equality. In~\cite{S20} the structure theorem was applied to the cyclic multiple planes $w^n=f(x,y)$ and to the varieties $\cV_m$ obtained as quotients of $m$-fold products of such a surface, with the aim of producing abelian varieties of large rank over the function fields $K_m=k(\cV_m)$. For curves this construction goes back to Yamagishi~\cite{Ya}; see \cite[Remark~5.1]{S19}.

The displayed statement is not correct, neither in its original form nor in the form of~\cite{S21}. A $K$-rational point of $\tilde A$ is a morphism $X'\to A$ that intertwines the covering automorphism with $\sigma$, and only the homomorphisms $\Prym_{X'/X}\to A$ with the same intertwining property contribute; the torsion part is also different. In the proof in \cite[\S4]{S19} the automorphism $\sigma$ does not enter the condition describing the points of the twist (Remark~\ref{rem:S19error}). Example~\ref{ex:counter}, with $A$ the elliptic curve $y^2=x^3+1$ over a number field containing $\zeta_3$, satisfies the hypotheses of \cite[Theorem~1.1]{S19} both before and after~\cite{S21}, and contradicts its conclusion. The error propagates to \cite[Theorem~1.2]{S19} (Example~\ref{ex:superelliptic}) and to the main results of~\cite{S20}: the ring $\End_k(\Alb(\cX_n))$ appearing there has to be replaced by the centralizer of the covering automorphism, several of the constants in \cite[Theorem~1.3]{S20} change, and some of the arguments need to be replaced.

The present paper replaces these statements by equivariant versions and gives complete proofs of the corrected results. Section~\ref{sec:corrections} records explicitly the changes to~\cite{S19}, \cite{S20} and~\cite{S21}. We also develop the arithmetic over number fields, which was only sketched in~\cite{S20}.

\subsection{Results}
Let $G$ be a finite group acting faithfully on a smooth, projective, geometrically integral variety $X'$ over a field $k$, put $L=k(X')$, $K=L^G$, and let $\rho\colon G\to\Aut_k(A)$ be an action of $G$ on an abelian variety $A$ over $k$ by automorphisms of abelian varieties. Let $\tilde A$ be the twist of $A$ by $L/K$ and $\rho$ (normalized as in Section~\ref{sec:twists}). Fix $x_0\in X'(k)$, let $\alpha\colon X'\to\Alb(X')$ be the Albanese map with $\alpha(x_0)=0$, and let $g\mapsto\tilde g$ be the induced action of $G$ on $\Alb(X')$.

\begin{thmx}[Theorem~\ref{thm:structure}]\label{thmA}
There is an exact sequence
\[
0\longrightarrow A(k)^G\longrightarrow \tilde A(K)\longrightarrow \Hom_k(\Alb(X'),A)^G\overset{\delta}{\longrightarrow} H^1(G,A(k)),
\]
where $\Hom_k(\Alb(X'),A)^G$ is the group of homomorphisms $\varphi$ with $\varphi\circ\tilde g=\rho(g)\circ\varphi$ for all $g\in G$. The image of $\delta$ is killed by $|G|$, and $\delta=0$ if $x_0$ is fixed by $G$, in which case the sequence splits. In particular
\begin{gather*}
\rk\tilde A(K)=\rk A(k)^G+\rk_\Z\Hom_k(\Alb(X'),A)^G,\\
\tilde A(K)_{\mathrm{tors}}=\bigl(A(k)^G\bigr)_{\mathrm{tors}}.
\end{gather*}
\end{thmx}

Now let $n\ge2$, let $k$ be a field of characteristic prime to $n$ containing a primitive $n$-th root of unity, and let $f\in k[x,y]$ be such that $w^n=f(x,y)$ defines a geometrically irreducible surface. Let $\cX_n$ be a smooth projective model, $\bar\tau$ the automorphism induced by $w\mapsto\zeta_nw$, $A=\Alb(\cX_n)$ and $\sigma=\bar\tau_*\in\Aut_k(A)$. For $m\ge1$ let $\cU_m=\cX_n^m$, let $\cV_m$ be the quotient of $\cU_m$ by the diagonal action of $\bar\tau$, and put $K_m=k(\cV_m)$. For an abelian variety $B$ over $k$ with an automorphism $\beta$ such that $\beta^n=\id$, let $\tilde B_m$ be the corresponding twist of $B$ over $K_m$; for $(B,\beta)=(A,\sigma)$ this is the Albanese variety of the twist $\tilde\cX_n$ considered in~\cite{S20}.

\begin{thmx}[Theorem~\ref{thm:planes}]\label{thmB}
Assume that $\cX_n(k)\ne\emptyset$. Then $1+\sigma+\dots+\sigma^{n-1}=0$, and
\begin{gather*}
\rk\tilde B_m(K_m)=\rk B^\beta(k)+m\cdot\rk_\Z\Hom_k(A,B)^{(\sigma,\beta)},\\
\tilde B_m(K_m)_{\mathrm{tors}}=B^\beta(k)_{\mathrm{tors}},
\end{gather*}
where $B^\beta=\ker(1-\beta)$ and $\Hom_k(A,B)^{(\sigma,\beta)}=\{\varphi:\varphi\circ\sigma=\beta\circ\varphi\}$. In particular
\[
\rk\Alb(\tilde\cX_n)(K_m)=m\cdot\rk_\Z\End_k(A)^{\sigma},\qquad \Alb(\tilde\cX_n)(K_m)_{\mathrm{tors}}=A^\sigma(k),
\]
where $\End_k(A)^\sigma$ is the centralizer of $\sigma$ in $\End_k(A)$. If $\cX_n$ has a $k$-rational point fixed by $\bar\tau$, then $\Alb(\tilde\cX_n)(K_m)$ is the direct sum of $A^\sigma(k)$ and a free left $\End_k(A)^\sigma$-module of rank $m$ with an explicit basis.
\end{thmx}

In~\cite{S20} the rank was claimed to be $m\cdot\rk\End_k(A)$. The two agree if and only if $\End_k(A)^\sigma$ has finite index in $\End_k(A)$; for instance, for $n=8$ and $n=12$ under the hypotheses of~\cite[Theorem~1.3]{S20}, this fails whenever $d_n>0$ and $k$ is large enough for all endomorphisms of $A_{\bar k}$ to be defined over $k$ (Remark~\ref{rem:n812}). Theorem~\ref{thmB} has the same form for cyclic covers $w^n=f(u_1,\dots,u_\ell)$ of $\PP^\ell$ whenever a smooth projective model carrying the induced cyclic action is available; in particular, this causes no difficulty in characteristic zero. For $\ell=1$ it gives the corrected form of \cite[Theorem~1.2]{S19} on superelliptic curves (Corollary~\ref{cor:superelliptic}), and Example~\ref{ex:superelliptic} shows that the torsion subgroup given there is not correct.

Over the complex numbers the centralizer is controlled by the eigenspace decomposition of $H_1(A,\Q)$ under $\sigma$ (Section~\ref{sec:hodge}). For $d\mid n$ let $e_d$ be the dimension over $\Q(\zeta_d)$ of the part of $H_1(A,\Q)$ on which $\sigma$ acts through primitive $d$-th roots of unity. The following theorem summarizes Corollary~\ref{cor:centralizer}, Theorem~\ref{thm:galois} and Proposition~\ref{prop:torsion}.

\begin{thmx}\label{thmC}
Let $k\subset\C$ be a number field containing $\zeta_n$, and assume $\cX_n(k)\ne\emptyset$.
\begin{enumerate}[\rm(i)]
\item $m\sum_{d\mid n,\ e_d>0}\phi(d)\ \le\ \rk\Alb(\tilde\cX_n)(K_m)\ \le\ m\sum_{d\mid n}\phi(d)\,e_d^2$.
\item The upper bound is attained if and only if $(A,\sigma)$ is rigid in the sense of Definition~\ref{def:rigid} and every endomorphism of $A_{\bar k}$ commuting with $\sigma$ is defined over $k$; the latter holds after a finite extension of $k$. Under Assumption~1.2 of~\cite{S20} the upper bound equals $m\cdot 4d_n^2/\phi(n)$, where $d_n=\dim A$.
\item $A^\sigma(\bar k)$ is a finite group of order $\prod p^{e_{p^r}}$, the product running over the prime powers $p^r>1$ dividing $n$. In particular $\Alb(\tilde\cX_n)(K_m)$ is torsion free if $e_d=0$ for every prime power $d$ dividing $n$.
\end{enumerate}
\end{thmx}

Part~(ii) replaces \cite[Theorem~1.3]{S20}: the constant there must be $4d_n^2/\phi(n)$ rather than $2d_n^2$ for $n=5,10$, and the equalities claimed there hold only when the relevant endomorphisms are defined over $k$. Section~\ref{sec:number} also contains a specialization statement (Proposition~\ref{prop:specialization}): for $m=1$ the base $\cV_1$ is rational, and N\'eron's specialization theorem produces, outside a thin set of $(a,b)\in k^2$, specializations attached to the $\mu_n$-torsors $w^n=f(a,b)$ whose Mordell--Weil rank over $k$ is at least $\rk\End_k(A)^\sigma$. For $m\ge2$ the corresponding assertion in \cite[Remark~1.5(2)]{S20} is not justified.

The dual of the Fermat cubic provides a concrete example in which the corrected theory can be worked out explicitly. Let
\[
F(x,y)=x^6+y^6+1-2(x^3y^3+x^3+y^3),
\]
so that $F=0$ is an affine equation of the dual of the Fermat cubic, a sextic with nine cusps.

\begin{thmx}[Theorem~\ref{thm:example}]\label{thmD}
\begin{enumerate}[\rm(i)]
\item The curve $Y^2=X^3+16F(x,y)$ has rank $3$ over $\Q(x,y)$ and rank $6$ over $\Q(\zeta_3)(x,y)$; the points $(4xy,4(x^3+y^3-1))$, $(4y,4(y^3+1-x^3))$, $(4x,4(x^3+1-y^3))$ generate a subgroup of finite index of its $\Q(x,y)$-rational points.
\item The curve $Y^2=X^3+F(x,y)$ has rank $0$ over $\Q(x,y)$ and over $\Q(\zeta_3)(x,y)$, rank $3$ over $\Q(\sqrt[3]{2})(x,y)$, and rank $6$ over $\Q(\zeta_3,\sqrt[3]{2})(x,y)$.
\item Let $A=\Alb(\cX_6)$ for $w^6=F(x,y)$. Then $\dim A=3$, Assumption~1.2 of~\cite{S20} holds, $A$ is isogenous over $\Q(\zeta_3)$ to the cube of $y^2=x^3+16$, and for every number field $k\supseteq\Q(\zeta_3)$ and every $m\ge1$ the group $\Alb(\tilde\cX_6)(k(\cV_m))$ is free of rank $18m$.
\item For all $(a,b)\in\Q^2$ outside a thin set, the Mordell curve $y^2=x^3+16F(a,b)$ has rank at least $3$ over $\Q$. The exceptional set is not empty.
\end{enumerate}
\end{thmx}

Part~(iii) shows that the rank formula of \cite[Theorem~1.4]{S20} does hold for the dual of the Fermat cubic, although the argument given there is incomplete. Part~(ii) illustrates the role of the field of definition: over $\Q(\zeta_3)$ the sextic $F$ has no quasi-toric decomposition $F=G^2+H^3$, whereas over $\Q(\sqrt[3]{2})$ it has three independent ones.

\subsection{Relation to other work}
Results relating Mordell--Weil groups of twists, or of isotrivial families, to groups of homomorphisms appear in several places; see for instance Hazama~\cite{Ha}, Wang~\cite{Wa}, Libgober~\cite{Li12}, and, for Jacobians over $k(t)$, Ulmer~\cite{Ul}. For $k=\C$ and $m=1$ the identification of the Mordell--Weil group with a group of equivariant morphisms is the starting point of the work of Cogolludo-Agust\'in and Libgober~\cite{CL} on the elliptic threefolds $u^2+v^3=F(x,y,1)$; Libgober~\cite{Li13} studies more generally isotrivial abelian varieties over $\C(x,y)$ with cyclic holonomy. In~\cite{CL} the rank is computed from the Alexander polynomial of the curve $F=0$ when this curve has only nodes and cusps. See also Hulek--Kloosterman~\cite{HK} and Kloosterman~\cite{Kl1,Kl2}. Theorem~\ref{thmA} is a version of this principle over an arbitrary field and for an arbitrary finite group, including the obstruction $\delta$ that appears when the cover has no rational fixed point. We include complete proofs because the statements in~\cite{S19} and~\cite{S20} have to be replaced. The observation that the growth of the rank in $m$ is formal (Proposition~\ref{prop:formal}) shows that the arithmetic content of the construction of~\cite{S20} lies in the case $m=1$; this is why the second half of the paper concentrates on the field of definition, torsion and specialization.

\subsection{Conventions}
Varieties are separated schemes of finite type over a field; points of a variety over a field $F$ are $F$-rational points. For an abelian variety $A$ over $k$ we write $\End_k(A)$ for the ring of endomorphisms defined over $k$. For any abelian group $M$ we use
\[
\rk M:=\dim_\Q(M\otimes_\Z\Q).
\]
This agrees with the usual rank when $M$ is finitely generated. For an automorphism $\sigma$ of $A$ we put $A^\sigma=\ker(1-\sigma)$, a closed subgroup scheme of $A$. We write $\phi$ for Euler's function and $\Phi_d$ for the $d$-th cyclotomic polynomial; the letter $\varphi$ is reserved for homomorphisms. For abelian varieties $A,B$ over $k$ with automorphisms $\sigma,\beta$,
\begin{gather*}
\Hom_k(A,B)^{(\sigma,\beta)}=\{\varphi\in\Hom_k(A,B):\varphi\circ\sigma=\beta\circ\varphi\},\\
\End_k(A)^\sigma=\Hom_k(A,A)^{(\sigma,\sigma)}.
\end{gather*}
We use the notation of~\cite{S20} whenever possible. References to numbered statements of~\cite{S19}, \cite{S20} and~\cite{S21} always carry the citation.

\section{Twists by Galois covers}\label{sec:twists}

\subsection{Conventions on twists}\label{subsec:conventions}
Let $L/K$ be a finite Galois extension and let $Y$ be a quasi-projective $K$-variety. An element $u\in\Gal(L/K)$ acts on $Y(L)$, and on morphisms defined over $L$, through its action on coordinates; we write ${}^uP$ and ${}^u\phi$. If $a=(a_u)$ is a $1$-cocycle of $\Gal(L/K)$ with values in $\Aut(Y_L)$ and $g\colon Y_L\to\tilde Y_L$ is the isomorphism of \cite[Prop.~3.1]{S20} (see \cite{BS}), so that ${}^ug=g\circ a_u$, then ${}^u(g(P))=g(a_u({}^uP))$ for $P\in Y(L)$, and therefore
\begin{equation}\label{eq:twistpoints}
\tilde Y(K)=\{\,g(P):P\in Y(L),\ a_u({}^uP)=P\ \text{for all }u\in\Gal(L/K)\,\}.
\end{equation}

\begin{setting}\label{set:cover}
Let $k$ be a field, $X'$ a smooth, projective, geometrically integral variety over $k$, and $G$ a finite group acting faithfully on $X'$ by $k$-automorphisms. Put $L=k(X')$ and $K=L^G$. For $g\in G$ let $u_g$ be the automorphism $h\mapsto h\circ g^{-1}$ of $L$. Since $u_gu_{g'}=u_{gg'}$ and $G$ acts faithfully on $L$, the map $g\mapsto u_g$ is an isomorphism $G\to\Gal(L/K)$. Let $A$ be an abelian variety over $k$ and $\rho\colon G\to\Aut_k(A)$ a homomorphism into the group of automorphisms of $A$ as an abelian variety. Since the elements of $\Aut_k(A)$ are fixed by $\Gal(L/K)$, the map $u_g\mapsto\rho(g)$ is a $1$-cocycle, and we denote by $\tilde A$ the corresponding twist of $A_K$; it exists because $A$ is projective, and it is an abelian variety over $K$.
\end{setting}

If $Y$ is a variety over $k$, the set $Y(L)$ is the set of rational maps $X'\dashrightarrow Y$ defined over $k$. If $f_P$ denotes the rational map corresponding to $P\in Y(L)$, then ${}^{u_g}P$ corresponds to $f_P\circ g^{-1}$, because the coordinates of $P$ are functions $h\in L$ and $u_g(h)=h\circ g^{-1}$. By~\eqref{eq:twistpoints}, if $Y$ carries an action $\rho$ of $G$ and $\tilde Y$ is the twist by $u_g\mapsto\rho(g)$, then
\begin{equation}\label{eq:equivariant}
\tilde Y(K)\ \longleftrightarrow\ \{\,f\colon X'\dashrightarrow Y\ \text{defined over }k : f\circ g=\rho(g)\circ f\ \text{for all }g\in G\,\}.
\end{equation}
For $Y=A$ this bijection is an isomorphism of groups, and every rational map $X'\dashrightarrow A$ is a morphism because $X'$ is smooth \cite[\S3]{Mi}. We always use this normalization; with the opposite identification of $G$ with $\Gal(L/K)$, $\rho$ has to be replaced by $g\mapsto\rho(g^{-1})$, which for cyclic $G$ amounts to replacing $\sigma$ by $\sigma^{-1}$.

Assume that $X'(k)\ne\emptyset$ and fix $x_0\in X'(k)$. Let $\alpha\colon X'\to\Alb(X')$ be the Albanese map with $\alpha(x_0)=0$ \cite{SeAlb}. By the universal property, every morphism $f\colon X'\to A$ defined over $k$ can be written uniquely as
\begin{equation}\label{eq:decomp}
f=c+\varphi\circ\alpha,\qquad c=f(x_0)\in A(k),\quad \varphi\in\Hom_k(\Alb(X'),A).
\end{equation}
For $g\in G$ let $\tilde g\in\Aut_k(\Alb(X'))$ be the automorphism defined by $\alpha\circ g=\tilde g\circ\alpha+\alpha(g(x_0))$; then $g\mapsto\tilde g$ is a homomorphism and does not depend on $x_0$. We write
\begin{gather*}
A(k)^G=\{c\in A(k):\rho(g)c=c\ \text{for all }g\in G\},\\
\Hom_k(\Alb(X'),A)^G=\{\varphi:\varphi\circ\tilde g=\rho(g)\circ\varphi\ \text{for all }g\in G\},
\end{gather*}
and we let $G$ act on $A(k)$ through $\rho$.

\subsection{The structure theorem}

\begin{theorem}\label{thm:structure}
In Setting~\ref{set:cover}, assume $X'(k)\ne\emptyset$ and fix $x_0\in X'(k)$.
\begin{enumerate}[\rm(a)]
\item There is an exact sequence
\[
0\longrightarrow A(k)^G\longrightarrow\tilde A(K)\overset{r}{\longrightarrow}\Hom_k(\Alb(X'),A)^G\overset{\delta}{\longrightarrow}H^1(G,A(k)),
\]
where $\delta(\varphi)$ is the class of the cocycle $g\mapsto-\varphi(\alpha(g(x_0)))$.
\item The image of $r$ contains $|G|\cdot\Hom_k(\Alb(X'),A)^G$. Consequently
\begin{gather*}
\rk\tilde A(K)=\rk A(k)^G+\rk_\Z\Hom_k(\Alb(X'),A)^G,\\
\tilde A(K)_{\mathrm{tors}}=\bigl(A(k)^G\bigr)_{\mathrm{tors}}.
\end{gather*}
\item If $x_0$ is fixed by $G$, then $\delta=0$ and
\[
\tilde A(K)\cong A(k)^G\oplus\Hom_k(\Alb(X'),A)^G.
\]
\end{enumerate}
\end{theorem}

\begin{proof}
Let $M=\Mor_k(X',A)=A(L)$, with $G$ acting by $g\cdot f=\rho(g)\circ f\circ g^{-1}$. By~\eqref{eq:equivariant}, $\tilde A(K)\cong M^G$. The constant maps form the $G$-submodule $A(k)$ of $M$, and by~\eqref{eq:decomp} the map $f\mapsto\varphi$ identifies $M/A(k)$ with $\Hom_k(\Alb(X'),A)$. This identification is $G$-equivariant for the action $g\cdot\varphi=\rho(g)\circ\varphi\circ\tilde g^{-1}$: indeed $\alpha\circ g^{-1}=\tilde g^{-1}\circ\alpha+\alpha(g^{-1}(x_0))$, so the homomorphism part of $\rho(g)\circ(c+\varphi\circ\alpha)\circ g^{-1}$ is $\rho(g)\circ\varphi\circ\tilde g^{-1}$. The invariants of $\Hom_k(\Alb(X'),A)$ under this action are exactly $\Hom_k(\Alb(X'),A)^G$. Taking $G$-invariants of
\[
0\to A(k)\to M\to\Hom_k(\Alb(X'),A)\to0
\]
gives the exact sequence in~(a). To compute $\delta$, lift $\varphi\in\Hom_k(\Alb(X'),A)^G$ to $f=\varphi\circ\alpha$. Then
\begin{align*}
g\cdot f-f&=\rho(g)\varphi\bigl(\tilde g^{-1}\alpha+\alpha(g^{-1}(x_0))\bigr)-\varphi\alpha\\
&=\rho(g)\varphi\bigl(\alpha(g^{-1}(x_0))\bigr)=\varphi\bigl(\tilde g\,\alpha(g^{-1}(x_0))\bigr),
\end{align*}
and $\tilde g\,\alpha(g^{-1}(x_0))=\alpha(x_0)-\alpha(g(x_0))=-\alpha(g(x_0))$. This proves~(a).

Since $H^1(G,A(k))$ is killed by $|G|$, we have $\delta(|G|\varphi)=0$, so $|G|\varphi$ lies in the image of $r$. As $\Hom_k(\Alb(X'),A)$ is a finitely generated free $\Z$-module, the rank formula follows, and since it is torsion free, the torsion of $\tilde A(K)$ is that of $A(k)^G$. This proves~(b). If $x_0$ is fixed by $G$, then $\alpha\circ g=\tilde g\circ\alpha$, so $\varphi\mapsto\varphi\circ\alpha$ is a $G$-equivariant section of $M\to\Hom_k(\Alb(X'),A)$, and~(c) follows.
\end{proof}

\begin{remark}\label{rem:prym}
Suppose that $N_\rho:=\sum_{g\in G}\rho(g)$ vanishes in $\End_k(A)$. If $\rho(g)c=c$ for all $g$ then $|G|c=N_\rho(c)=0$, so $A^G:=\bigcap_g\ker(1-\rho(g))$ is contained in $A[|G|]$ and $A(k)^G$ is finite. Also, every $\varphi\in\Hom_k(\Alb(X'),A)^G$ satisfies $\varphi\circ\sum_g\tilde g=N_\rho\circ\varphi=0$, hence factors through the quotient $\Prym_{X'/X}=\Alb(X')/\operatorname{Im}(\sum_g\tilde g)$, which for cyclic $G$ is the Prym variety of \cite[Def.~2.2]{S20}. Thus
\[
\rk\tilde A(K)=\rk_\Z\Hom_k(\Prym_{X'/X},A)^G.
\]
\end{remark}

For cyclic groups we obtain the following statement, which replaces \cite[Theorem~1.1]{S19}, its amended form in \cite[(2)]{S21}, and \cite[Theorem~3.3]{S20}.

\begin{corollary}\label{cor:cyclic}
In Setting~\ref{set:cover} let $G=\langle\gamma\rangle$ be cyclic of order $n$, put $\sigma=\rho(\gamma)$ and $\tilde\gamma$ as above, and assume $X'(k)\ne\emptyset$.
\begin{enumerate}[\rm(a)]
\item $\rk\tilde A(K)=\rk A^\sigma(k)+\rk_\Z\Hom_k(\Alb(X'),A)^{(\tilde\gamma,\sigma)}$ and $\tilde A(K)_{\mathrm{tors}}=A^\sigma(k)_{\mathrm{tors}}$. If $\gamma$ has a $k$-rational fixed point, then $\tilde A(K)\cong A^\sigma(k)\oplus\Hom_k(\Alb(X'),A)^{(\tilde\gamma,\sigma)}$.
\item If $1+\sigma+\dots+\sigma^{n-1}=0$, then $A^\sigma\subset A[n]$ and $\rk\tilde A(K)=\rk_\Z\Hom_k(\Prym_{X'/X},A)^{(\tilde\gamma,\sigma)}$.
\item If in addition there are an abelian variety $B$ over $k$ with an automorphism $\beta$, and a $k$-isogeny $\lambda\colon\Prym_{X'/X}\to A^m\times B$ with $\lambda\circ\tilde\gamma=(\sigma\times\dots\times\sigma\times\beta)\circ\lambda$ and $\Hom_k(B,A)^{(\beta,\sigma)}=0$, then $\rk\tilde A(K)=m\cdot\rk_\Z\End_k(A)^\sigma$.
\end{enumerate}
\end{corollary}

\begin{proof}
Parts (a) and (b) are Theorem~\ref{thm:structure} and Remark~\ref{rem:prym}. For~(c), put $\theta=\sigma\times\dots\times\sigma\times\beta$ and let $\lambda'$ be a $k$-isogeny with $\lambda'\lambda=[N]$ and $\lambda\lambda'=[N]$. From $\lambda\tilde\gamma=\theta\lambda$ we get $\lambda'\theta\lambda=N\tilde\gamma=\tilde\gamma\lambda'\lambda$, hence $\lambda'\theta=\tilde\gamma\lambda'$ because $\lambda$ is surjective. Composition with $\lambda'$ and with $\lambda$ therefore gives homomorphisms between $\Hom_k(\Prym_{X'/X},A)^{(\tilde\gamma,\sigma)}$ and $\Hom_k(A^m\times B,A)^{(\theta,\sigma)}$ whose composites are multiplication by $N$, and the latter group is $(\End_k(A)^\sigma)^m\oplus\Hom_k(B,A)^{(\beta,\sigma)}=(\End_k(A)^\sigma)^m$.
\end{proof}

\begin{remark}\label{rem:S19error}
In \cite[\S4]{S19} a point of $\tilde A(K)$ is written as a pair $(\lambda,Q)\in\Hom_k(\Alb(X'),A)\oplus A(k)$, as in~\eqref{eq:decomp}, on which $\gamma^j$ is said to act by $(\lambda,Q)\mapsto(\lambda\circ\tilde\gamma^j,Q)$, and the condition for $(\lambda,Q)$ to give a point of the twist is written as $\lambda\circ\tilde\gamma^j=\lambda\circ\tilde\gamma^{n-j}$ for all $j$; this is then replaced by the conditions that $\lambda$ vanish on the image of $1+\tilde\gamma+\dots+\tilde\gamma^{n-1}$ and that $Q\in A[n](k)$. The automorphism $\sigma$ defining the twist does not enter these conditions, and the formula for the action ignores the translation by $\lambda(\alpha(\gamma^j(x_0)))$, which vanishes when the base point is fixed by $\gamma$ but not in general. By~\eqref{eq:equivariant} and~\eqref{eq:decomp}, the correct conditions are $\lambda\circ\tilde\gamma=\sigma\circ\lambda$ and, when the base point is fixed by $\gamma$, $\sigma(Q)=Q$. In particular the answer depends on $\sigma$, whereas the statement of \cite[Theorem~1.1]{S19} does not.
\end{remark}

\begin{remark}\label{rem:n2}
For $n=2$ and $\sigma=-1$ the automorphism induced by $\tilde\gamma$ on $\Prym_{X'/X}$ is $-1$, so every homomorphism $\Prym_{X'/X}\to A$ is equivariant and $A^\sigma=A[2]$. In this case, which is the setting of quadratic twists considered by Hazama and recalled in \cite[Theorem~3.1]{S19}, Corollary~\ref{cor:cyclic} agrees with the rank statements of \cite[Theorem~3.1]{S19} and \cite[Theorem~3.3]{S20}, and gives the isomorphism $\tilde A(K)\cong\Hom_k(\Prym_{X'/X},A)\oplus A[2](k)$ when $\gamma$ has a $k$-rational fixed point. For $n\ge3$ the statements differ in general, as the following example shows.
\end{remark}

\subsection{A counterexample}

\begin{example}\label{ex:counter}
Let $k\subset\C$ be a field containing a primitive cube root of unity $\rho_3$. Let $X'=E$ be the elliptic curve $y^2=x^3+1$ with origin $O$ at infinity, and $\gamma(x,y)=(\rho_3x,y)$. Then $\gamma$ is an automorphism of order $3$ fixing $O$, hence a group automorphism, and $\End_k(E)=\Z[\gamma]\cong\Z[\rho_3]$ is commutative and has no zero divisors; in particular $1+\gamma+\gamma^2=0$. The function $t=y$ generates the function field of $E/\langle\gamma\rangle\cong\PP^1$, so $K=k(t)$ and $L=K(x)$ with $x^3=t^2-1$. Here $\Alb(E)=E$, $\alpha=\id$, $\tilde\gamma=\gamma$ and $\Prym_{E/\PP^1}=E$. Take $A=E$ and $\sigma\in\{\gamma,\gamma^2\}$. The hypotheses of \cite[Theorem~3.3]{S20} hold with $m=1$ and $B=0$, and that theorem asserts that $\rk\tilde A(K)=\rk\End_k(E)=2$ for both choices of $\sigma$. If $k$ is a number field, the hypotheses of \cite[Theorem~1.1]{S19} also hold (with $s=3$, $n=1$ and $\dim B=0$ in the notation there), so that \cite{S19} asserts $\rk\tilde A(K)\ge2$ and \cite[(2)]{S21} asserts $\rk\tilde A(K)=2$, for both choices of $\sigma$.

By Corollary~\ref{cor:cyclic}(a), applied with the fixed point $O$, we have $\tilde A(K)\cong E^\sigma(k)\oplus\End_k(E)^{(\gamma,\sigma)}$. For $\sigma=\gamma$ this gives rank $2$. For $\sigma=\gamma^2$ the condition $\varphi\gamma=\gamma^2\varphi$ reads $(\gamma-\gamma^2)\varphi=0$, whence $\varphi=0$ and $\tilde A(K)=E^\sigma(k)$ is finite.

Explicitly, the twist is $Y^2=X^3+(t^2-1)^2$ for $\sigma=\gamma$, with $g(x_0,y_0)=(x^2x_0,x^3y_0)$, and $Y^2=X^3+(t^2-1)^4$ for $\sigma=\gamma^2$, with $g(x_0,y_0)=(x^4x_0,x^6y_0)$; one checks ${}^{u_\gamma}g=g\circ\sigma$ using $u_\gamma(x)=\rho_3^{-1}x$ and $u_\gamma(t)=t$. The vanishing of the rank of the second curve can also be seen geometrically: over $\C(t)$ its Weierstrass equation is minimal and defines an elliptic K3 surface with singular fibres of type $\mathrm{IV}^*$ at $t=\pm1$ and $t=\infty$ \cite{Ta}; the trivial lattice has rank $2+3\cdot6=20$, the maximal Picard number of a complex K3 surface, so the Mordell--Weil rank is $0$ by the Shioda--Tate formula \cite{Sh}. The first curve defines a rational elliptic surface with three fibres of type $\mathrm{IV}$, of Mordell--Weil rank $10-2-6=2$.

Finally $E^\sigma(k)=\{O,(0,1),(0,-1)\}\cong\Z/3\Z$. If $\sqrt[3]{2}\in k$, then $E[3](k)\cong(\Z/3\Z)^2$, the $3$-division polynomial of $E$ being $3x(x^3+4)$. For $\sigma=\gamma^2$ we then have $\tilde A(K)\cong\Z/3\Z$, while \cite[Theorem~1.1]{S19} and \cite[Theorem~3.3]{S20} predict $\Z^2\oplus(\Z/3\Z)^2$; for $\sigma=\gamma$ we have $\tilde A(K)\cong\Z^2\oplus\Z/3\Z$.
\end{example}

\section{Cyclic multiple planes}\label{sec:planes}

\subsection{Setting}\label{subsec:planes-setting}
Let $n\ge2$ and let $k$ be a field of characteristic prime to $n$ containing a primitive $n$-th root of unity $\zeta_n$. Let $f\in k[x,y]$ be a polynomial such that the affine surface $X_n\colon w^n=f(x,y)$ is geometrically irreducible, let $\cX_n$ be a smooth projective model of $X_n$ on which the automorphism $\tau(x,y,w)=(x,y,\zeta_nw)$ induces an automorphism $\bar\tau$, as in \cite[\S4]{S20}, and put
\[
A=\Alb(\cX_n),\qquad \sigma=\bar\tau_*\in\Aut_k(A),\qquad d_n=\dim A,
\]
where $\bar\tau_*$ is the automorphism induced by $\bar\tau$ (Section~\ref{subsec:conventions}, with $G=\langle\bar\tau\rangle$). Any two smooth projective models of $X_n$ are birational, so $A$ and $\sigma$ do not depend on the choice of $\cX_n$ up to isomorphism, and by the Lang--Nishimura lemma \cite{Ni} neither does the condition $\cX_n(k)\ne\emptyset$; it holds, for instance, if $X_n$ has a non-singular $k$-rational point. For $m\ge1$ let $\cU_m=\cX_n^{(1)}\times_k\dots\times_k\cX_n^{(m)}$ with coordinates $(x_i,y_i,w_i)$ on the $i$-th factor, let $\bar\gamma=(\bar\tau,\dots,\bar\tau)$, let $\cV_m=\cU_m/\langle\bar\gamma\rangle$, and put $L_m=k(\cU_m)$ and $K_m=k(\cV_m)=L_m^{\langle\bar\gamma\rangle}$. Let $B$ be an abelian variety over $k$ and $\beta\in\Aut_k(B)$ with $\beta^n=\id$. We denote by $\tilde B_m$ the twist of $B$ over $K_m$ defined, as in Setting~\ref{set:cover}, by $G=\langle\bar\gamma\rangle$ acting on $\cU_m$ and $\rho(\bar\gamma)=\beta$. Similarly $\tilde\cX_n$ denotes the twist of $\cX_n$ defined by $\rho(\bar\gamma)=\bar\tau$; by \cite[Prop.~3.2]{S20}, $\Alb(\tilde\cX_n)$ is the twist $\tilde A_m$ of $A$ defined by $\rho(\bar\gamma)=\sigma$.

\begin{lemma}\label{lem:twistmodel}
The twist $\tilde\cX_n$ is birational over $K_m$ to the surface $f(x_1,y_1)z^n=f(x,y)$ of \cite[(5.3)]{S20}. Under~\eqref{eq:equivariant} the $K_m$-rational points $P_1=(x_1,y_1,1)$ and $P_i=(x_i,y_i,w_i/w_1)$ $(2\le i\le m)$ of \cite[(5.4)]{S20} correspond to the projections $\pr_i\colon\cU_m\to\cX_n^{(i)}=\cX_n$.
\end{lemma}

\begin{proof}
Let $u=u_{\bar\gamma}$. Since $w_1\circ\bar\gamma^{-1}=\zeta_n^{-1}w_1$, the $L_m$-isomorphism $g(x,y,w)=(x,y,w/w_1)$ from $X_n$ to $f(x_1,y_1)z^n=f(x,y)$ satisfies ${}^ug(x,y,w)=(x,y,\zeta_nw/w_1)=g(\tau(x,y,w))$, that is ${}^ug=g\circ\tau$. The point $P_i$ is the image under $g$ of the $L_m$-point $(x_i,y_i,w_i)$ of $X_n$, which is the rational map $\pr_i$.
\end{proof}

\begin{lemma}\label{lem:norm}
We have $1+\sigma+\dots+\sigma^{n-1}=0$ in $\End_k(A)$. Consequently $\Prym_{\cX_n/\PP^2}=A$, and $A^\sigma\subset A[n]$ is finite.
\end{lemma}

\begin{proof}
It suffices to prove the identity over $\bar k$. Fix $x_0\in\cX_n(\bar k)$ and let $\alpha$ be the Albanese map with $\alpha(x_0)=0$. The morphism $\Phi=\sum_{j=0}^{n-1}\alpha\circ\bar\tau^j$ satisfies $\Phi\circ\bar\tau=\Phi$. Regarded as a point of $A(\bar k(\cX_n))$, it is therefore invariant under $\Gal(\bar k(\cX_n)/\bar k(x,y))$, which is generated by $h\mapsto h\circ\bar\tau^{-1}$ because $X_n$ is geometrically irreducible of degree $n$ over the plane. Hence $\Phi$ factors through a rational map $\PP^2\dashrightarrow A$, which is constant \cite[\S3]{Mi}. On the other hand $\Phi=(1+\sigma+\dots+\sigma^{n-1})\circ\alpha+c$ for some $c\in A(\bar k)$, and $\alpha(\cX_n)$ generates $A$. The last assertions follow from Remark~\ref{rem:prym}.
\end{proof}

\begin{remark}
In \cite[\S5]{S20} this identity was deduced from $(1-\sigma)(1+\sigma+\dots+\sigma^{n-1})=0$ and $\sigma\ne1$. That deduction is not valid, since $\End(A)$ may have zero divisors (for instance when $A$ is isogenous to $E^{d}$ with $d\ge2$).
\end{remark}

\subsection{Mordell--Weil groups over $K_m$}

\begin{theorem}\label{thm:planes}
In the setting of~\S\ref{subsec:planes-setting}, assume $\cX_n(k)\ne\emptyset$.
\begin{enumerate}[\rm(a)]
\item There is an exact sequence
\[
0\to B^\beta(k)\to\tilde B_m(K_m)\to\bigl(\Hom_k(A,B)^{(\sigma,\beta)}\bigr)^m\to H^1(\Z/n\Z,B(k)),
\]
and the image of the middle map contains $n\cdot(\Hom_k(A,B)^{(\sigma,\beta)})^m$. Hence
\begin{gather*}
\rk\tilde B_m(K_m)=\rk B^\beta(k)+m\cdot\rk_\Z\Hom_k(A,B)^{(\sigma,\beta)},\\
\tilde B_m(K_m)_{\mathrm{tors}}=B^\beta(k)_{\mathrm{tors}}.
\end{gather*}
\item For $(B,\beta)=(A,\sigma)$, the group $A^\sigma(k)$ is finite and
\[
\rk\Alb(\tilde\cX_n)(K_m)=m\cdot\rk_\Z\End_k(A)^\sigma,\qquad \Alb(\tilde\cX_n)(K_m)_{\mathrm{tors}}=A^\sigma(k).
\]
\item Assume that $\cX_n$ has a $k$-rational point $O$ fixed by $\bar\tau$ (for example the point $(a,b,0)$, where $(a,b)$ is a non-singular $k$-rational point of the curve $f=0$, if $\cX_n$ is chosen so that $\cX_n\dashrightarrow X_n$ is an isomorphism over the smooth locus of $X_n$). Then $O$ defines a $K_m$-rational point $\tilde O$ of $\tilde\cX_n$; let $\tilde Q_i$ be the image of $P_i$ under the Albanese map of $\tilde\cX_n$ with base point $\tilde O$. Every $\varphi\in\Hom_k(A,B)^{(\sigma,\beta)}$ descends to a $K_m$-homomorphism $\tilde\varphi\colon\Alb(\tilde\cX_n)\to\tilde B_m$, and
\[
\tilde B_m(K_m)=B^\beta(k)\oplus\bigoplus_{i=1}^m\{\tilde\varphi(\tilde Q_i):\varphi\in\Hom_k(A,B)^{(\sigma,\beta)}\},
\]
where each summand $\varphi\mapsto\tilde\varphi(\tilde Q_i)$ is injective. In particular $\Alb(\tilde\cX_n)(K_m)=A^\sigma(k)\oplus\bigoplus_{i=1}^m\End_k(A)^\sigma\cdot\tilde Q_i$, and $\tilde Q_1,\dots,\tilde Q_m$ is a basis of a free left $\End_k(A)^\sigma$-module.
\end{enumerate}
\end{theorem}

\begin{proof}
Apply Theorem~\ref{thm:structure} to $X'=\cU_m$, $G=\langle\bar\gamma\rangle\cong\Z/n\Z$ and $\rho(\bar\gamma)=\beta$, with base point $(x_0,\dots,x_0)$ for some $x_0\in\cX_n(k)$ (and $x_0=O$ in~(c)). Since $\Alb(\cX_n\times\cX_n)=A\times A$ with Albanese map $\alpha\times\alpha$ \cite{SeAlb}, we have $\Alb(\cU_m)=A^m$, the Albanese map is $\alpha^{\times m}$, and $\tilde{\bar\gamma}=\sigma^{\times m}$. Therefore
\[
\Hom_k(A^m,B)^{(\sigma^{\times m},\beta)}=\bigoplus_{i=1}^m\Hom_k(A,B)^{(\sigma,\beta)}\circ e_i ,
\]
where $e_i\colon A^m\to A$ is the $i$-th projection. This gives~(a). For~(b), $A^\sigma(k)$ is finite by Lemma~\ref{lem:norm}.

In~(c), $\bar\tau(O)=O$ gives $\alpha\circ\bar\tau=\sigma\circ\alpha$, so the constant map with value $O$ corresponds to a point $\tilde O\in\tilde\cX_n(K_m)$ by~\eqref{eq:equivariant}, and $\alpha$ descends to the Albanese map $\tilde\alpha$ of $\tilde\cX_n$ with base point $\tilde O$ (see the proof of \cite[Prop.~3.2]{S20}). By Lemma~\ref{lem:twistmodel}, $\tilde Q_i=\tilde\alpha(P_i)$ corresponds to $\alpha\circ\pr_i=e_i\circ\alpha^{\times m}$, that is, to the pair $(c,\varphi)=(0,e_i)$ in~\eqref{eq:decomp}. If $\varphi\in\Hom_k(A,B)^{(\sigma,\beta)}$, then $g_B\circ\varphi\circ g_A^{-1}$ is invariant under $\Gal(L_m/K_m)$, where $g_A$ and $g_B$ are the isomorphisms of~\eqref{eq:twistpoints}; this is the descended homomorphism $\tilde\varphi$, and it acts on points by $f\mapsto\varphi\circ f$ in terms of~\eqref{eq:equivariant}. Thus $\tilde\varphi(\tilde Q_i)$ corresponds to $(0,\varphi\circ e_i)$, and~(c) follows from Theorem~\ref{thm:structure}(c).
\end{proof}

\begin{remark}\label{rem:rankformula}
Theorem~\ref{thm:planes}(b) replaces \cite[Theorem~1.1]{S20}, which asserts $\rk\Alb(\tilde\cX_n)(K_m)=m\cdot\rk\End_k(A)$ and has $\Alb(\cX_n)[m](k)$ as torsion part. The two rank formulas agree if and only if $\End_k(A)^\sigma$ has finite index in $\End_k(A)$, for instance if $\sigma$ is central in $\End_k(A)\otimes\Q$. The assertion at the end of \cite[\S5]{S20} that $\tilde Q_1,\dots,\tilde Q_m$ form a set of independent generators is replaced by Theorem~\ref{thm:planes}(c).
\end{remark}

The dependence on $m$ is formal, as the following description of $K_m$ shows. Let $\tilde\cX_n^{(1)}$ denote the twist of $\cX_n$ over $K_1=k(x_1,y_1)$, which by Lemma~\ref{lem:twistmodel} is birational to $f(x_1,y_1)z^n=f(x,y)$.

\begin{proposition}\label{prop:formal}
The field $K_m$ is $K_1$-isomorphic to the function field of $(\tilde\cX_n^{(1)})^{m-1}$ over $K_1$. Consequently, if $\cX_n(k)\ne\emptyset$ and $(B,\beta)=(A,\sigma)$, the rank over $K_m$ is $m$ times the rank over $K_1$.
\end{proposition}

\begin{proof}
Put $u_i=w_i/w_1\in L_m$ for $2\le i\le m$; these are invariant under $\bar\gamma$, and $f(x_1,y_1)u_i^n=f(x_i,y_i)$. Let $Y$ be the $(m-1)$-fold product of the surface $f(x_1,y_1)z^n=f(x,y)$ over $K_1$; it is geometrically integral since $X_n$ is. The point $(x_i,y_i,u_i)_{2\le i\le m}$ of $Y$ has transcendence degree $2(m-1)$ over $K_1$, hence is a generic point, so the subfield $F=k(x_1,y_1,x_i,y_i,u_i:2\le i\le m)$ of $K_m$ is $K_1$-isomorphic to $K_1(Y)$. Since $L_m=F(w_1)$ with $w_1^n=f(x_1,y_1)\in F$, we have $[L_m:F]\le n=[L_m:K_m]$, hence $F=K_m$. The last assertion is Theorem~\ref{thm:planes}(b).
\end{proof}

Proposition~\ref{prop:formal} is an instance of the general fact that $B'(F(Y))\cong B'(F)\oplus\Hom_F(\Alb(Y),B')$ for an abelian variety $B'$ over a field $F$ and a smooth, projective, geometrically integral variety $Y$ over $F$ with an $F$-point: over $F=K_1$, the variety $Y=(\tilde\cX_n^{(1)})^{m-1}$ contributes $m-1$ further copies of $\End_{K_1}(\tilde A_1)$. The arithmetic content of the construction is therefore concentrated in the case $m=1$, that is, in the twist of $A$ over $k(x,y)$.

\subsection{Cyclic covers of $\PP^\ell$ and superelliptic curves}\label{subsec:otherdims}
Let $\ell\ge1$ and let $f\in k[u_1,\dots,u_\ell]$ be such that the hypersurface $w^n=f(u_1,\dots,u_\ell)$ is geometrically irreducible. Suppose that it admits a smooth projective model on which $w\mapsto\zeta_n w$ extends. Replacing $(x,y)$ by $(u_1,\dots,u_\ell)$, the arguments of this section then apply without change: the only property of $\PP^2$ used in Lemma~\ref{lem:norm} is that the quotient is rational. The same applies to Sections~\ref{sec:hodge} and~\ref{sec:number}; in Proposition~\ref{prop:specialization} one uses that $\AAA^\ell(k)$ is not thin. The existence of such a smooth equivariant model is automatic in characteristic zero. For $\ell\ge3$ these are the covers considered in \cite[Remark~1.5(1)]{S20}.

\begin{remark}\label{rem:higherdim}
In particular \cite[Remark~1.5(1)]{S20} holds in the form of Theorem~\ref{thm:planes}.
\end{remark}

For $\ell=1$ we recover the setting of \cite[\S5]{S19}. Let $f\in k[x]$ be a polynomial of degree $r\ge1$ with non-zero discriminant, so that $f$ is square-free and the curve $y^n=f(x)$ is geometrically irreducible. Let $C$ be its smooth projective model, $\iota$ the automorphism induced by $y\mapsto\zeta_ny$, $J=J(C)=\Alb(C)$, and $\tilde\iota\in\Aut_k(J)$ the induced automorphism. For $m\ge1$ let $\cV_m=C^m/\langle(\iota,\dots,\iota)\rangle$ and $K_m=k(\cV_m)$, and let $C^\xi$ be the twist of $C$ over $K_m$ considered in \cite[\S5]{S19}. By Lemma~\ref{lem:twistmodel}, $C^\xi$ is birational to $f(x_1)z^n=f(x)$, which is the model given in \cite[\S5]{S19}, so that $C^\xi$ is the twist of $C$ by $\rho(\bar\gamma)=\iota$ and $J(C^\xi)$ is the twist of $J$ by $\rho(\bar\gamma)=\tilde\iota$. Let $P_1,\dots,P_m$ be the $K_m$-rational points of \cite[\S5, (1)]{S19} (with the indexing corrected as in Section~\ref{sec:corrections}); they are the points of Lemma~\ref{lem:twistmodel}.

\begin{corollary}[Corrected form of {\cite[Theorem~1.2]{S19}}]\label{cor:superelliptic}
Assume that $C(k)\ne\emptyset$. Then $J^{\tilde\iota}\subset J[n]$ is finite,
\[
\rk J(C^\xi)(K_m)=m\cdot\rk_\Z\End_k(J)^{\tilde\iota},\qquad J(C^\xi)(K_m)_{\mathrm{tors}}=J^{\tilde\iota}(k).
\]
If $C$ has a $k$-rational point $O$ fixed by $\iota$, for instance $O=(a,0)$ with $a\in k$ a root of $f$, then
\[
J(C^\xi)(K_m)=J^{\tilde\iota}(k)\oplus\bigoplus_{i=1}^m\End_k(J)^{\tilde\iota}\cdot\tilde Q_i,
\]
where $\tilde Q_i$ is the image of $P_i$ under the Abel--Jacobi map of $C^\xi$ with base point $\tilde O$. In particular, the rank formula $\rk J(C^\xi)(K_m)=m\cdot\rk\End_k(J)$ of \cite[(3)]{S21} holds if and only if $\tilde\iota$ is central in $\End_k(J)\otimes\Q$; this is the case for $n=2$ and whenever $\End_k(J)\otimes\Q$ is commutative.
\end{corollary}

\begin{proof}
This is Theorem~\ref{thm:planes}(b),(c) and Lemma~\ref{lem:norm} for $\ell=1$ and $(B,\beta)=(J,\tilde\iota)$; the point $(a,0)$ is a non-singular point of the affine curve because $f$ is square-free. The rank of $\End_k(J)^{\tilde\iota}$ equals that of $\End_k(J)$ if and only if the centralizer of $\tilde\iota$ in $\End_k(J)\otimes\Q$ is everything, that is, if and only if $\tilde\iota$ is central.
\end{proof}

In \cite[\S5]{S19} the identity $1+\tilde\iota+\dots+\tilde\iota^{n-1}=0$ is deduced from $(1-\tilde\iota)(1+\tilde\iota+\dots+\tilde\iota^{n-1})=0$ and $\tilde\iota\neq1$; as in~\cite{S20}, this deduction is not valid in general, but the identity holds by Lemma~\ref{lem:norm}, since $C/\langle\iota\rangle\cong\PP^1$. The assertion at the end of \cite[\S5]{S19} that $\tilde Q_1,\dots,\tilde Q_m$ belong to a set of independent generators is replaced by Corollary~\ref{cor:superelliptic}.

\begin{example}\label{ex:superelliptic}
Let $n=3$, $f=x^2-1$ and $m\ge2$, and let $k$ be a number field containing $\zeta_3$ and $\sqrt[3]2$. These data satisfy the hypotheses of \cite[Theorem~1.2]{S19}, which requires $\deg f\le m$, and $c=(1,0)$ is a $k$-rational point of $C$ fixed by $\iota$. The map $(x,y)\mapsto(y,x)$ is an isomorphism from $C\colon y^3=x^2-1$ to the curve $E\colon y^2=x^3+1$ of Example~\ref{ex:counter}, carrying $\iota$ to $\gamma$ (with $\rho_3=\zeta_3$). Hence $J\cong E$, $\End_k(J)\cong\Z[\zeta_3]$ is commutative, $J^{\tilde\iota}(k)\cong E^\gamma(k)\cong\Z/3\Z$, and $J[3](k)\cong(\Z/3\Z)^2$ by Example~\ref{ex:counter}. By Corollary~\ref{cor:superelliptic},
\[
J(C^\xi)(K_m)\cong\Z^{2m}\oplus\Z/3\Z,
\]
whereas \cite[Theorem~1.2]{S19} asserts $J(C^\xi)(K_m)\cong\Z^{2m}\oplus(\Z/3\Z)^2$. For $m=1$, which is excluded in \cite{S19}, the model $(x_1^2-1)z^3=x^2-1$ of $C^\xi$ is isomorphic over $K_1=k(x_1)$ to $Y^2=X^3+(x_1^2-1)^2$ via $X=(x_1^2-1)z$, $Y=(x_1^2-1)x$; this is the twist of Example~\ref{ex:counter} for $\sigma=\gamma$, in accordance with Lemma~\ref{lem:twistmodel}.
\end{example}

\section{Endomorphisms commuting with an automorphism}\label{sec:hodge}

In this section all abelian varieties are over $\C$. We describe groups of equivariant homomorphisms in terms of the eigenspace decomposition of the first homology. This replaces the case-by-case analysis of \cite[\S6]{S20}.

\subsection{Eigenspaces}\label{subsec:eigen}
Let $A$ be a complex abelian variety and $\sigma\in\Aut(A)$ with $\sigma^n=\id$. Put $\Lambda_A=H_1(A,\Z)$ and $V_A=\Lambda_A\otimes\Q$. The tangent space $T_0A$ is $V_A\otimes\R$ with a complex structure $J$; let $W_A\subset V_A\otimes\C$ be the $i$-eigenspace of $J\otimes1$. Then $V_A\otimes\C=W_A\oplus\overline{W_A}$, and $v\mapsto\frac12(v-iJv)$ is a $\sigma$-equivariant $\C$-linear isomorphism $T_0A\to W_A$. In the notation of \cite[\S1]{S20}, $\Lambda_n=\Lambda_A$ and $W_n=W_A$. Recall that a $\Q$-linear map $u\colon V_A\to V_B$ comes from $\Hom(A,B)\otimes\Q$ if and only if $u_\C(W_A)\subseteq W_B$.

For $d\mid n$ let $V_{A,d}=\ker\Phi_d(\sigma)$, where $\Phi_d$ is the $d$-th cyclotomic polynomial. Since $\sigma$ has finite order, $V_A=\bigoplus_{d\mid n}V_{A,d}$. Fix a primitive $n$-th root of unity $\zeta_n\in\C$ and put $\zeta_d=\zeta_n^{n/d}$. Letting $\zeta_d$ act as $\sigma$ makes $V_{A,d}$ a vector space over $\Q(\zeta_d)$; we put $e_d(A)=\dim_{\Q(\zeta_d)}V_{A,d}$, so that $2\dim A=\sum_{d\mid n}\phi(d)e_d(A)$. For $\lambda\in\C$ with $\lambda^n=1$, of order $d$, let $\iota_\lambda\colon\Q(\zeta_d)\to\C$ be the embedding with $\iota_\lambda(\zeta_d)=\lambda$, and let
\begin{gather*}
H_A(\lambda)=V_{A,d}\otimes_{\Q(\zeta_d),\iota_\lambda}\C,\\
W_A(\lambda)=W_A\cap H_A(\lambda),\qquad a_A(\lambda)=\dim_\C W_A(\lambda).
\end{gather*}
Then $H_A(\lambda)$ is the $\lambda$-eigenspace of $\sigma$ on $V_A\otimes\C$ and has dimension $e_d(A)$. Since $W_A$ and $\overline{W_A}$ are stable under $\sigma$, we have $W_A=\bigoplus_\lambda W_A(\lambda)$ and $H_A(\lambda)=W_A(\lambda)\oplus(\overline{W_A}\cap H_A(\lambda))$; complex conjugation maps $W_A(\bar\lambda)$ onto $\overline{W_A}\cap H_A(\lambda)$. Hence
\begin{equation}\label{eq:conj}
a_A(\lambda)+a_A(\bar\lambda)=e_d(A).
\end{equation}

\begin{definition}\label{def:rigid}
We say that $(A,\sigma)$ is \emph{rigid} if $a_A(\lambda)\in\{0,e_d(A)\}$ for every $\lambda$, where $d$ is the order of $\lambda$. In that case, for every $d$ with $e_d(A)>0$ we put $\Sigma_A(d)=\{\lambda:\ \lambda\text{ of order }d,\ W_A(\lambda)=H_A(\lambda)\}$; by~\eqref{eq:conj}, the primitive $d$-th roots of unity are the disjoint union of $\Sigma_A(d)$ and $\overline{\Sigma_A(d)}$.
\end{definition}

If $\lambda=\pm1$ then $\lambda=\bar\lambda$ and~\eqref{eq:conj} gives $a_A(\lambda)=e_d(A)/2$; hence rigidity forces $e_1(A)=e_2(A)=0$. Assumption~1.2 of \cite{S20} requires that $\Lambda_n$ be a $\Z[\zeta_n]$-module, that is $V_A=V_{A,n}$, and that no pair of complex conjugate characters occur in $W_n$, that is, $W_A(\lambda)$ and $W_A(\bar\lambda)$ are not both non-zero. By~\eqref{eq:conj}, Assumption~1.2 holds if and only if $(A,\sigma)$ is rigid and $V_A=V_{A,n}$.

\subsection{Equivariant homomorphisms}

\begin{lemma}\label{lem:equivhom}
Let $(A,\sigma)$ and $(B,\beta)$ be complex abelian varieties with automorphisms satisfying $\sigma^n=\beta^n=\id$. Then $\Hom(A,B)^{(\sigma,\beta)}\otimes\Q$ is the set of
\[
u\in\bigoplus_{d\mid n}\Hom_{\Q(\zeta_d)}(V_{A,d},V_{B,d})
\]
such that $u_\lambda(W_A(\lambda))\subseteq W_B(\lambda)$ for all $\lambda$, where $u_\lambda=u\otimes_{\iota_\lambda}\C\colon H_A(\lambda)\to H_B(\lambda)$.
\end{lemma}

\begin{proof}
A $\Q$-linear map $u\colon V_A\to V_B$ satisfies $u\sigma=\beta u$ if and only if it maps each $V_{A,d}$ into $V_{B,d}$ and is $\Q(\zeta_d)$-linear there. Such a map sends $H_A(\lambda)$ into $H_B(\lambda)$, and since $W_B=\bigoplus_\lambda W_B(\lambda)$, the condition $u_\C(W_A)\subseteq W_B$ is equivalent to $u_\lambda(W_A(\lambda))\subseteq W_B(\lambda)$ for all $\lambda$.
\end{proof}

\begin{corollary}\label{cor:equivhom}
In the situation of Lemma~\ref{lem:equivhom}:
\begin{enumerate}[\rm(i)]
\item $\dim_\Q\Hom(A,B)^{(\sigma,\beta)}\otimes\Q\le\sum_{d\mid n}\phi(d)\,e_d(A)\,e_d(B)$, with equality if and only if, for every $\lambda$, either $W_A(\lambda)=0$ or $W_B(\lambda)=H_B(\lambda)$.
\item If $(A,\sigma)$ and $(B,\beta)$ are rigid, then
\[
\dim_\Q\Hom(A,B)^{(\sigma,\beta)}\otimes\Q=\sum_{d\in D}\phi(d)\,e_d(A)\,e_d(B),
\]
where $D$ is the set of $d\mid n$ with $e_d(A)e_d(B)>0$ and $\Sigma_A(d)=\Sigma_B(d)$.
\end{enumerate}
\end{corollary}

\begin{proof}
The space
\[
S=\bigoplus_d\Hom_{\Q(\zeta_d)}(V_{A,d},V_{B,d})
\]
has dimension $\sum_d\phi(d)e_d(A)e_d(B)$, and $u\mapsto(u_\lambda)_\lambda$ identifies $S\otimes\C$ with $\bigoplus_\lambda\Hom_\C(H_A(\lambda),H_B(\lambda))$. The conditions of Lemma~\ref{lem:equivhom} are $\C$-linear in $(u_\lambda)$, so they hold on all of $S$ if and only if they hold on all of $S\otimes\C$, that is, if and only if every linear map $H_A(\lambda)\to H_B(\lambda)$ sends $W_A(\lambda)$ into $W_B(\lambda)$. This proves~(i).

For~(ii), the conditions decouple according to $d$. Fix $d$ with $e_d(A)e_d(B)>0$. If $\Sigma_A(d)=\Sigma_B(d)$, then for $\lambda\in\Sigma_A(d)$ we have $W_B(\lambda)=H_B(\lambda)$, and for $\lambda\notin\Sigma_A(d)$ we have $W_A(\lambda)=0$; so the whole $d$-component of $S$ satisfies the conditions. Otherwise, since $\Sigma_A(d)$ and $\Sigma_B(d)$ both have $\phi(d)/2$ elements, there is $\lambda\in\Sigma_A(d)\smallsetminus\Sigma_B(d)$; then $W_A(\lambda)=H_A(\lambda)$ and $W_B(\lambda)=0$, so $u_\lambda=0$, and this forces $u|_{V_{A,d}}=0$ because $u\mapsto u\otimes_{\iota_\lambda}\C$ is injective.
\end{proof}

\begin{corollary}\label{cor:centralizer}
Let $(A,\sigma)$ be a complex abelian variety with an automorphism, $\sigma^n=\id$. Then the centralizer of $\sigma$ in $\End(A)\otimes\Q$ contains $\Q[\sigma]\cong\prod_{e_d(A)>0}\Q(\zeta_d)$ and is contained in $\bigoplus_d\End_{\Q(\zeta_d)}(V_{A,d})\cong\prod_d M_{e_d(A)}(\Q(\zeta_d))$. Hence
\[
\sum_{d\mid n,\ e_d(A)>0}\phi(d)\ \le\ \rk_\Z\End(A)^\sigma\ \le\ \sum_{d\mid n}\phi(d)\,e_d(A)^2,
\]
and the upper bound is attained if and only if $(A,\sigma)$ is rigid. In particular, if Assumption~1.2 of \cite{S20} holds and $n\ge3$, then the centralizer is isomorphic to $M_e(\Q(\zeta_n))$ with $e=2\dim A/\phi(n)$, and
\[
\rk_\Z\End(A)^\sigma=\frac{4(\dim A)^2}{\phi(n)}.
\]
\end{corollary}

\begin{proof}
Apply Corollary~\ref{cor:equivhom}(i) with $(B,\beta)=(A,\sigma)$: the equality condition says that $W_A(\lambda)\in\{0,H_A(\lambda)\}$ for all $\lambda$, which is rigidity. The lower bound holds because $\sigma$ is semisimple and its minimal polynomial is $\prod_{e_d(A)>0}\Phi_d$. The last statement is the case $V_A=V_{A,n}$.
\end{proof}

For $n=3,4,6$ this gives $2(\dim A)^2$, and for $n=5,8,10,12$ it gives $(\dim A)^2$. Note that the restriction $3\le n\le 12$, $n\ne7,9,11$ of \cite[Theorem~1.3]{S20} plays no role.

\subsection{The cases $n=8$ and $n=12$}
The following lemma is elementary; we include a proof because it is needed to correct \cite[Prop.~4.2(v),(vi)]{S20}.

\begin{lemma}\label{lem:quadratic}
Let $B$ be a complex abelian variety of dimension $g$ and $F_0\subset\End(B)\otimes\Q$ an imaginary quadratic field acting on $T_0B$ through a single embedding $\iota_0\colon F_0\to\C$. Then $B$ is isogenous to $E_0^g$, where $E_0$ is an elliptic curve with complex multiplication by the ring of integers $\cO_{F_0}$ of $F_0$, and $\End(B)\otimes\Q\cong M_g(F_0)$.
\end{lemma}

\begin{proof}
The space $V=H_1(B,\Q)$ is an $F_0$-vector space of dimension $g$. Since $F_0$ acts on $T_0B=V\otimes_\Q\R$ by the scalars $\iota_0(a)$, the complex structure of $V\otimes_\Q\R$ is multiplication by the element of $F_0\otimes_\Q\R$ mapped to $i$ by $\iota_0\otimes\R$. An $F_0$-basis of $V$ therefore gives a $\C$-linear isomorphism $T_0B\cong\C^g$ under which $H_1(B,\Z)$ is commensurable with $\iota_0(\cO_{F_0})^g$. Hence $B$ is isogenous to $(\C/\iota_0(\cO_{F_0}))^g=E_0^g$, and $\End(B)\otimes\Q\cong M_g(F_0)$.
\end{proof}

\begin{remark}\label{rem:n812}
Let $n\in\{8,12\}$ and let $(A,\sigma)$ satisfy Assumption~1.2 of \cite{S20}. By rigidity and $V_A=V_{A,n}$, we have $W_A=\bigoplus_{\lambda\in\Sigma_A(n)}H_A(\lambda)$ and $\Sigma_A(n)=\{\zeta_n^j:j\in S\}$ for a subset $S\subset(\Z/n\Z)^\times$ with $(\Z/n\Z)^\times=S\sqcup(-S)$. For $n=8$ the possible sets $S$ are $\{1,3\}$, $\{5,7\}$, $\{1,5\}$, $\{3,7\}$, i.e.\ the cosets of $\{1,3\}$ and $\{1,5\}$, which are the subgroups fixing $\Q(\zeta_8+\zeta_8^3)=\Q(\sqrt{-2})$ and $\Q(\zeta_8^2)=\Q(i)$. For $n=12$ they are the cosets $\{1,7\},\{5,11\}$ of the subgroup fixing $\Q(\zeta_{12}^4)=\Q(\sqrt{-3})$ and $\{1,5\},\{7,11\}$ of the subgroup fixing $\Q(\zeta_{12}^3)=\Q(i)$. Hence there is an imaginary quadratic field $F_0\subset\Q(\zeta_n)=\Q[\sigma]$ such that all the embeddings $\iota_\lambda$, $\lambda\in\Sigma_A(n)$, agree on $F_0$; thus $F_0$ acts on $T_0A\cong W_A$ through a single embedding, and Lemma~\ref{lem:quadratic} shows that $A$ is isogenous to $E_0^{\dim A}$ with $E_0$ having complex multiplication by $\cO_{F_0}$. Consequently:
\begin{enumerate}[\rm(1)]
\item $\End(A)\otimes\Q\cong M_{\dim A}(F_0)$ has dimension $2(\dim A)^2$, whereas the centralizer of $\sigma$ has dimension $(\dim A)^2$ by Corollary~\ref{cor:centralizer}. So, if $\dim A>0$, then $\sigma$ is not central, and $m\cdot\rk\End(A)$ is not the rank.
\item In \cite[Prop.~4.2(v),(vi)]{S20} one of $n_1,n_2$ vanishes. For~(vi) this follows from the above, since $E_i$ and $E_\rho$ are not isogenous. For~(v), let $\zeta_8$ be a square root of $i$ and $s(x,y)=(ix,\zeta_8y)$. This is an automorphism of order $8$ of $C_2\colon y^2=x^5+x$, and $s^4$ is the hyperelliptic involution, so $s_*^4=-1$ on $J(C_2)$ and $\Q[s_*]\cong\Q(\zeta_8)$. The automorphism $s$ acts on the basis $dx/y$, $x\,dx/y$ of $H^0(C_2,\Omega^1)$ by $\zeta_8$ and $\zeta_8^3$, hence $s_*$ acts on $T_0J(C_2)=H^0(C_2,\Omega^1)^\vee$ with the same eigenvalues, and the corresponding embeddings of $\Q(\zeta_8)$ agree on $\Q(\zeta_8+\zeta_8^3)=\Q(\sqrt{-2})$. By Lemma~\ref{lem:quadratic}, $J(C_2)$ is isogenous to $E_{\sqrt{-2}}^2$, where $E_{\sqrt{-2}}$ has complex multiplication by $\Z[\sqrt{-2}]$ and is not isogenous to $E_i$.
\end{enumerate}
\end{remark}

\begin{remark}\label{rem:alexander}
For a cyclic multiple plane, $e_1(A)=0$ by Lemma~\ref{lem:norm}. The numbers $e_d(A)$ and $a_A(\lambda)$ can be computed from the Alexander polynomial of the branch curve and from the position of its singular points; see Libgober~\cite{Li82}, Esnault~\cite{Es} and, for curves whose only singularities are nodes and cusps, Cogolludo-Agust\'in and Libgober~\cite{CL}. We use this only through~\cite{CL} in Section~\ref{sec:example}.
\end{remark}

\section{Arithmetic over number fields}\label{sec:number}

We return to the setting of~\S\ref{subsec:planes-setting}, now with $k$ of characteristic $0$. Let $\bar k$ be an algebraic closure and $\Gamma_k=\Gal(\bar k/k)$. For a finite extension $k'$ of $k$ in $\bar k$ we write $K'_m=k'(\cV_m)$.

\subsection{The Galois action}

\begin{theorem}\label{thm:galois}
Assume $\cX_n(k)\ne\emptyset$, and let $H=\Hom_{\bar k}(A,B)^{(\sigma,\beta)}$. Let $k_H$ be the fixed field of the kernel of the natural action $\Gamma_k\to\Aut(H)$. Then $k_H/k$ is finite, the action factors through $\Gal(k_H/k)$, and for every finite extension $k'$ of $k$ in $\bar k$
\[
\rk\tilde B_m(K'_m)=\rk B^\beta(k')+m\cdot\rk_\Z H^{\Gamma_{k'}}.
\]
Moreover $\rk_\Z H^{\Gamma_{k'}}=\rk_\Z H$ if and only if $k'\supseteq k_H$.
\end{theorem}

\begin{proof}
Theorem~\ref{thm:planes}(a) over $k'$ gives the formula with $\Hom_{k'}(A,B)^{(\sigma,\beta)}$ in place of $H^{\Gamma_{k'}}$; these groups coincide because $\sigma$ and $\beta$ are defined over $k$. The group $H$ is finitely generated and each of its elements is defined over a finite extension of $k$, so $\Gamma_k$ acts through a finite quotient. Finally, if $H^{\Gamma_{k'}}$ has finite index in $H$, then every element of $\Gamma_{k'}$ acts trivially on $H\otimes\Q$, hence on the torsion free group $H$.
\end{proof}

\begin{corollary}\label{cor:number}
Let $k\subset\C$ be a number field containing $\zeta_n$, assume $\cX_n(k)\neq\emptyset$, and let $e_d=e_d(A_\C)$. Then
\[
m\sum_{d\mid n,\ e_d>0}\phi(d)\ \le\ \rk\Alb(\tilde\cX_n)(K_m)\ \le\ m\sum_{d\mid n}\phi(d)\,e_d^2.
\]
The upper bound is attained if and only if $(A_\C,\sigma)$ is rigid and every endomorphism of $A_{\bar k}$ commuting with $\sigma$ is defined over $k$. If Assumption~1.2 of \cite{S20} holds and $n\ge3$, then
\[
\rk\Alb(\tilde\cX_n)(K_m)\le m\cdot\frac{4d_n^2}{\phi(n)},
\]
with equality if and only if $k\supseteq k_H$ for $H=\End_{\bar k}(A)^\sigma$; if moreover $2d_n=\phi(n)$, equality holds over $k$.
\end{corollary}

\begin{proof}
Fix an embedding $\bar k\hookrightarrow\C$ extending the given embedding $k\hookrightarrow\C$. Endomorphisms of abelian varieties do not change after extension between algebraically closed fields, so $\End_{\bar k}(A)=\End(A_\C)$. The lower bound holds because $\Z[\sigma]\subseteq\End_k(A)^\sigma$, and the rest follows from Theorem~\ref{thm:planes}(b), Theorem~\ref{thm:galois} and Corollary~\ref{cor:centralizer}. If $2d_n=\phi(n)$ under Assumption~1.2, then $e=1$ and the centralizer is $\Q(\zeta_n)=\Q[\sigma]$, which is defined over $k$.
\end{proof}

Rigidity of $(A_\C,\sigma)$ does not depend on the embedding $k\subset\C$: the numbers $a_A(\lambda)$ are the dimensions of the eigenspaces of $\sigma$ on the $k$-vector space $\mathrm{Lie}(A)$, whose eigenvalues lie in $k$, and a change of embedding only permutes the labels $\lambda$. Corollary~\ref{cor:number} replaces \cite[Theorem~1.3]{S20}. It also gives the correct form of \cite[Theorem~1.4]{S20}: under the hypotheses there, $\alpha(\cX_6)$ is a surface by \cite[Theorem~4.3]{S20}, which is a result of Kulikov \cite{Ku}, so $d_6\ge2$, and $\rk\Alb(\tilde\cX_6)(K_m)\le 2md_6^2$, with equality if and only if every endomorphism of $A_{\bar k}$ commuting with $\sigma$ is defined over $k$.

\subsection{Torsion}

\begin{proposition}\label{prop:torsion}
Let $k\subset\C$ and assume $\cX_n(k)\ne\emptyset$. Then $\Alb(\tilde\cX_n)(K_m)_{\mathrm{tors}}=A^\sigma(k)$, and $A^\sigma(\bar k)\cong\Lambda_A/(1-\sigma)\Lambda_A$ is a group killed by $n$, of order
\[
\bigl|A^\sigma(\bar k)\bigr|=\prod_{p^r\mid n,\ r\ge1}p^{\,e_{p^r}},
\]
the product running over the prime powers $p^r>1$ dividing $n$. In particular $\Alb(\tilde\cX_n)(K_m)$ is torsion free if $e_d=0$ for every prime power $d>1$ dividing $n$; this is the case, for example, if Assumption~1.2 of \cite{S20} holds and $n$ is not a prime power.
\end{proposition}

\begin{proof}
The first statement is Theorem~\ref{thm:planes}(b). By Lemma~\ref{lem:norm}, $e_1=0$, so $1-\sigma$ is invertible on $V_A$. Writing $A(\C)=(V_A\otimes\R)/\Lambda_A$, the kernel of $1-\sigma$ is $(1-\sigma)^{-1}\Lambda_A/\Lambda_A$, which $1-\sigma$ maps isomorphically onto $\Lambda_A/(1-\sigma)\Lambda_A$; since $A^\sigma$ is finite, $A^\sigma(\bar k)=A^\sigma(\C)$. Its order is $|\det(1-\sigma\mid V_A)|=\prod_{d}\Phi_d(1)^{e_d}$, and $\Phi_d(1)$ equals $p$ if $d=p^r$ with $r\ge1$ and equals $1$ if $d>1$ is not a prime power. Finally $n=\sum_{j=1}^{n-1}(1-\sigma^j)$ in $\End(A)$ by Lemma~\ref{lem:norm}, so $n$ kills $\Lambda_A/(1-\sigma)\Lambda_A$.
\end{proof}

\begin{remark}\label{rem:cokernel}
If $\cX_n$ has no $k$-rational point fixed by $\bar\tau$, the map $\delta$ of Theorem~\ref{thm:structure} takes values in $H^1(\Z/n\Z,A(k))$. Since $1+\sigma+\dots+\sigma^{n-1}=0$, one has $H^1(\Z/n\Z,A(k))=A(k)/(1-\sigma)A(k)$, a finite group when $k$ is a number field, because it is killed by $n$ and $A(k)$ is finitely generated. The image of $\delta$ is an arithmetic invariant; it vanishes over an algebraically closed field, where the isogeny $1-\sigma$ is surjective on points.
\end{remark}

\subsection{Specialization}
Let $k$ be a number field containing $\zeta_n$, let $U\subset\AAA^2_k$ be the open set $f\ne0$, and let $T\subset U\times\AAA^1$ be the $\mu_n$-torsor $w^n=f(x,y)$, with $\mu_n$ acting by $w\mapsto\zeta w$. For $c\in k^\times$ let $T_c=\Spec k[w]/(w^n-c)$ and let $B_c$ be the twist of $B$ by $T_c$ and $\beta$: its points over an extension $k'$ of $k$ are the $\mu_n$-equivariant morphisms $T_c\otimes k'\to B$, where $\zeta_n$ acts on $B$ by $\beta$. Up to isomorphism $B_c$ depends only on the class of $c$ in $k^\times/k^{\times n}$.

\begin{proposition}\label{prop:specialization}
Assume $\cX_n(k)\ne\emptyset$.
\begin{enumerate}[\rm(a)]
\item There is an abelian scheme $\cB\to U$ whose generic fibre is $\tilde B_1$ and whose fibre at $u\in U(k)$ is $B_{f(u)}$.
\item There is a thin set $\Omega\subset U(k)$ such that for $u\in U(k)\smallsetminus\Omega$ the specialization map $\tilde B_1(K_1)\to B_{f(u)}(k)$ is injective. In particular
\[
\rk B_{f(u)}(k)\ \ge\ \rk B^\beta(k)+\rk_\Z\Hom_k(A,B)^{(\sigma,\beta)}\qquad(u\in U(k)\smallsetminus\Omega).
\]
\item The classes of $f(u)$ in $k^\times/k^{\times n}$, for $u\in U(k)\smallsetminus\Omega$, form an infinite set.
\end{enumerate}
\end{proposition}

\begin{proof}
(a) The torsor $T\to U$ is finite \'etale because $n$ and $f$ are invertible on $U$. Choose an ample line bundle $\mathcal L$ on $B$ and replace it by $\bigotimes_{j=0}^{n-1}(\beta^j)^*\mathcal L$, which is ample and $\beta$-invariant. Descent of $B\times T$ along $T\to U$ with respect to the descent datum defined by $\beta$ is then effective \cite[\S6.1]{BLR}, and gives an abelian scheme $\cB\to U$. Its points over a $U$-scheme $S$ are the $\mu_n$-equivariant $S$-morphisms $T\times_US\to B\times S$, so its generic fibre is $\tilde B_1$ by~\eqref{eq:equivariant} and Lemma~\ref{lem:twistmodel}, and its fibre at $u$ is $B_{f(u)}$ because $T_u=T_{f(u)}$.

(b) This is N\'eron's specialization theorem \cite{Ne}, in the form of \cite[\S11.1]{SeMW}, applied to $\cB\to U$; the rank of $\tilde B_1(K_1)$ is given by Theorem~\ref{thm:planes}(a) with $m=1$.

(c) For $c\in k^\times$, the set $\Omega_c$ of $u\in U(k)$ with $f(u)\in c\,k^{\times n}$ is the image of the $k$-points of the surface $cz^n=f(x,y)$, which is geometrically irreducible and of degree $n\ge2$ over $U$; so $\Omega_c$ is thin \cite[\S3.1]{SeTopics}. If only finitely many classes $c_1,\dots,c_r$ occurred, then $U(k)$ would be contained in the thin set $\Omega\cup\Omega_{c_1}\cup\dots\cup\Omega_{c_r}$, which contradicts Hilbert's irreducibility theorem \cite[Chap.~3]{SeTopics}.
\end{proof}

\begin{remark}\label{rem:mge2}
For $m\ge2$ the rational points of $\cV_m$ lie, on a dense open set, over $m$-tuples of points $(x_i,y_i)\in k^2$ such that $f(x_i,y_i)/f(x_1,y_1)\in k^{\times n}$, and the specialization of $\tilde B_m$ at such a point is again $B_{f(x_1,y_1)}$. N\'eron's theorem gives injectivity of specialization outside a thin subset of $\cV_m(k)$, which is useful only if $\cV_m(k)$ is not thin. We do not know whether this holds in general, and it was not established in~\cite{S20}; accordingly \cite[Remark~1.5(2)]{S20} should be restricted to the case $m=1$, which is Proposition~\ref{prop:specialization}. Silverman's theorem \cite{Si}, quoted there, concerns families over a curve.
\end{remark}

\subsection{Elliptic targets with $j=0$}
Let $n=6$, write $\zeta=\zeta_6$, and for $D\in k^\times$ let $E_D$ be the elliptic curve $y^2=x^3+D$ with the automorphism $\beta(x,y)=(\zeta^{-2}x,-y)$ of order $6$. Since $\beta^*(dx/y)=\zeta\,dx/y$, the automorphism $\beta$ acts on $T_0E_D$ by $\zeta$; thus $(E_D,\beta)$ is rigid with $\Sigma_{E_D}(6)=\{\zeta\}$. The fixed points of $\beta$ are trivial, so $E_D^\beta(k)=0$. The $L_1$-isomorphism $g(x_0,y_0)=(w^2x_0,w^3y_0)$ from $E_D$ to
\[
Y^2=X^3+Df(x,y)
\]
satisfies ${}^{u}g=g\circ\beta$ for $u=u_{\bar\tau}$, since $u(w)=\zeta^{-1}w$. Hence $Y^2=X^3+Df(x,y)$ is the twist $\tilde B_1$ for $(B,\beta)=(E_D,\beta)$, and Theorem~\ref{thm:planes} gives, if $\cX_6(k)\ne\emptyset$,
\begin{equation}\label{eq:mordellrank}
\rk E_{Df}(K_1)=\rk_\Z\Hom_k(A,E_D)^{(\sigma,\beta)},\qquad E_{Df}(K_1)_{\mathrm{tors}}=0.
\end{equation}
A decomposition $Df=G^2+H^3$ with $G,H\in k(x,y)$ gives the point $(-H,G)$. Over $K_m$, with $u_i=w_i/w_1$, the same decomposition gives the points $\bigl(-H(x_i,y_i)/u_i^2,\ G(x_i,y_i)/u_i^3\bigr)$ of $Y^2=X^3+Df(x_1,y_1)$, since $u_i^6=f(x_i,y_i)/f(x_1,y_1)$.

\begin{lemma}\label{lem:threeisog}
Let $F_1$ be a field of characteristic $0$ with $\sqrt{-3}\notin F_1$, and $c\in F_1^\times$. Then $\dim_\Q E_c(F_1(\sqrt{-3}))\otimes\Q=2\dim_\Q E_c(F_1)\otimes\Q$.
\end{lemma}

\begin{proof}
For the quadratic extension $F_1(\sqrt{-3})/F_1$ one has $E_c(F_1(\sqrt{-3}))\otimes\Q\cong E_c(F_1)\otimes\Q\oplus E_c^{(-3)}(F_1)\otimes\Q$, where $E_c^{(-3)}\colon y^2=x^3-27c$ is the quadratic twist by $-3$. The map $(x,y)\mapsto\bigl((x^3+4c)/x^2,\ y(x^3-8c)/x^3\bigr)$ is an isogeny $E_c\to E_{-27c}$ of degree $3$ defined over $F_1$, with kernel $\{O,(0,\pm\sqrt c)\}$, so the two summands have the same dimension.
\end{proof}

\section{The dual of the Fermat cubic}\label{sec:example}

In this section $n=6$, $\zeta=\zeta_6$, $\rho_3=\zeta^2$, and
\begin{gather*}
F(x,y)=x^6+y^6+1-2(x^3y^3+x^3+y^3),\\
\bar F(u,v,w)=u^6+v^6+w^6-2(u^3v^3+v^3w^3+w^3u^3).
\end{gather*}
The Gauss map of the Fermat cubic $x^3+y^3+z^3=0$ sends $(x:y:z)$ to the dual point $(x^2:y^2:z^2)$. If $a=x^3$, $b=y^3$ and $c=z^3$, then
\[
\bar F(x^2,y^2,z^2)=(a-b-c)(a-b+c)(a+b-c)(a+b+c),
\]
so the image of the Gauss map lies on $\bar F=0$. Since the dual of a smooth plane cubic is an irreducible sextic, $C\colon\bar F=0$ is the dual of the Fermat cubic, and $F(x,y)=\bar F(x,y,1)$.

The singular locus can also be read off directly. The partial derivatives are
\[
\bar F_u=6u^2(u^3-v^3-w^3),\qquad
\bar F_v=6v^2(v^3-w^3-u^3),\qquad
\bar F_w=6w^2(w^3-u^3-v^3).
\]
These equations show that a singular point has one zero coordinate and that the cubes of the other two coordinates are equal. Thus the singular points are exactly
\[
(0:1:\epsilon),\qquad (1:0:\epsilon),\qquad (1:\epsilon:0),
\qquad \epsilon^3=1.
\]
At $(0:1:1)$, on the chart $v=1$ with $s=w-1$, the local equation has the form
\[
9s^2-4u^3+R(u,s),\qquad R(u,s)\in (s^3,su^3,u^6).
\]
Thus the tangent cone is the double line $s=0$, while the restriction to this tangent has order $3$; the singularity is an ordinary cusp. The group generated by permutations of the coordinates and by multiplication of a coordinate by a cube root of unity preserves $\bar F$ and acts transitively on the nine singular points, so all nine are ordinary cusps. Its Alexander polynomial is $(t^2-t+1)^3$; see \cite[Introduction]{CL} and \cite[Theorem~1.2.5]{De}.

Put $g_1=x^3+y^3-1$, $g_2=y^3+1-x^3$, $g_3=x^3+1-y^3$ and $h_1=xy$, $h_2=y$, $h_3=x$. Then
\begin{equation}\label{eq:decomps}
F=g_i^2-4h_i^3\qquad(i=1,2,3).
\end{equation}
Consequently the points
\[
Q_i=(4h_i,\,4g_i)\qquad(i=1,2,3)
\]
lie on $E_{16F}\colon Y^2=X^3+16F$ over $\Q(x,y)$, and the points $P_i=(\sqrt[3]{4}\,h_i,\,g_i)$ lie on $E_F\colon Y^2=X^3+F$ over $\Q(\sqrt[3]2)(x,y)$. The map $(X,Y)\mapsto(\sqrt[3]{16}\,X,4Y)$ is an isomorphism $E_F\to E_{16F}$ over $\Q(\sqrt[3]2)(x,y)$ sending $P_i$ to $Q_i$. We write $[\rho_3]$ for the automorphism $(X,Y)\mapsto(\rho_3X,Y)$ of these curves.

\begin{theorem}\label{thm:example}
\begin{enumerate}[\rm(i)]
\item $E_{16F}(\Q(x,y))$ has rank $3$, and $Q_1,Q_2,Q_3$ generate a subgroup of finite index. $E_{16F}(\Q(\zeta_3)(x,y))$ has rank $6$, and $Q_i,[\rho_3]Q_i$ $(i=1,2,3)$ generate a subgroup of finite index.
\item $E_F$ has rank $0$ over $\Q(x,y)$ and over $\Q(\zeta_3)(x,y)$, rank $3$ over $\Q(\sqrt[3]2)(x,y)$ and rank $6$ over $\Q(\zeta_3,\sqrt[3]2)(x,y)$. Moreover $E_F(\Q(\zeta_3)(x,y))=0$; in particular $F$ has no decomposition $F=G^2+H^3$ with $G,H\in\Q(\zeta_3)(x,y)$.
\item Let $A=\Alb(\cX_6)$ for $w^6=F(x,y)$ and $\sigma=\bar\tau_*$. Then $\dim A=3$, $(A,\sigma)$ satisfies Assumption~1.2 of \cite{S20}, and $A$ is isogenous over $\Q(\zeta_3)$, compatibly with $\sigma$ and $\beta\times\beta\times\beta$, to $E_{16}^3$, where $E_{16}\colon y^2=x^3+16$. For every number field $k\supseteq\Q(\zeta_3)$ and every $m\ge1$, the group $\Alb(\tilde\cX_6)(k(\cV_m))$ is free of rank $18m=2md_6^2$.
\item For all $(a,b)\in\Q^2$ outside a thin set, $F(a,b)\ne0$ and the points $Q_i(a,b)$ are independent in the group of rational points of the Mordell curve $y^2=x^3+16F(a,b)$, which therefore has rank at least $3$ over $\Q$. For $(a,b)=(2,3)$ one has $Q_1+Q_2=2Q_3$, so the exceptional set is not empty.
\end{enumerate}
\end{theorem}

\begin{proof}
\emph{Geometric rank.} By \cite[Theorem~1.1]{CL}, the rank of $u^2+v^3=F(x,y)$ over $\C(x,y)$ equals $\deg\Delta_C=6$. (The statement there is for an affine equation $\bar F(x,y,1)$; changing the line at infinity multiplies $\bar F$ by the sixth power of a rational function, which does not change the curve over $\C(\PP^2)$ up to isomorphism.) Since $16$ is a sixth power in $\C$, the curves $E_F$ and $E_{16F}$ are isomorphic over $\C(x,y)$, so both have rank $6$ over $\C(x,y)$ and at most $6$ over any subfield.

\emph{Independence.} The Weierstrass equation of $E_{16F}$ defines an elliptic curve over the open set $U\colon F\ne0$ of $\AAA^2_\Q$, and $Q_1,Q_2,Q_3$ are sections over $U$, since their coordinates are polynomials. Because $U$ is integral, a relation $\sum n_iQ_i=0$ over $\Q(x,y)$ holds for the sections and hence in every fibre. Take $(a,b)=(2,5)$, so that $16F(2,5)=214784=2^8\cdot839$ and
\[
Q_1(2,5)=(40,528),\qquad Q_2(2,5)=(20,472),\qquad Q_3(2,5)=(8,-464)
\]
on $E\colon y^2=x^3+214784$. For $p\in\{19,41,67\}$ the curve $E$ has good reduction, with
\[
|E(\F_{19})|=28,\qquad |E(\F_{41})|=42,\qquad |E(\F_{67})|=63.
\]
Thus the $7$-primary part is cyclic of order $7$ in each case. Take generators
\[
R_{19}=(14,15),\qquad R_{41}=(6,23),\qquad R_{67}=(40,8)
\]
of these order-$7$ subgroups. After multiplication by the prime-to-$7$ cofactors $4,6,9$, respectively, the reductions of $Q_1(2,5),Q_2(2,5),Q_3(2,5)$ have coefficient matrix
\[
\begin{pmatrix}
1&6&5\\
1&3&6\\
5&4&4
\end{pmatrix}
\quad\text{over }\F_7,
\]
where the rows correspond to $p=19,41,67$. Its determinant is $5$ modulo $7$, so the three classes are independent in the product of the three quotients. Since $|E(\F_5)|=6$, good reduction at $5$ also gives $E(\Q)[7]=0$. If $\sum n_iQ_i(2,5)=0$ with not all $n_i$ zero, we may divide the coefficients by their common power of $7$; the absence of rational $7$-torsion shows that the resulting relation is still zero. Reducing modulo the three primes then contradicts the nonsingularity of the matrix above. Hence $Q_1(2,5),Q_2(2,5),Q_3(2,5)$ are independent, and so are $Q_1,Q_2,Q_3$.

\emph{Proof of~(i).} Let $c$ denote the non-trivial automorphism of the extension $\Q(\zeta_3)(x,y)/\Q(x,y)$. The $\Q$-vector space $\mathcal M=E_{16F}(\Q(\zeta_3)(x,y))\otimes\Q$ is the direct sum of its $(\pm1)$-eigenspaces for $c$, and the $(+1)$-eigenspace is $E_{16F}(\Q(x,y))\otimes\Q$. The endomorphism $[\sqrt{-3}]=[\rho_3]-[\rho_3^2]$ is an isogeny satisfying $c\circ[\sqrt{-3}]=-[\sqrt{-3}]\circ c$, so it maps the $(+1)$-eigenspace injectively into the $(-1)$-eigenspace. Hence $Q_1,Q_2,Q_3,[\sqrt{-3}]Q_1,[\sqrt{-3}]Q_2,[\sqrt{-3}]Q_3$ are independent, and so are the $Q_i,[\rho_3]Q_i$, since $[\rho_3]=(-1+[\sqrt{-3}])/2$. Together with the upper bound $6$, this gives rank $6$ over $\Q(\zeta_3)(x,y)$, with the $Q_i,[\rho_3]Q_i$ spanning $\mathcal M$, and rank $3$ over $\Q(x,y)$ with the $Q_i$ spanning the $(+1)$-eigenspace, in accordance with Lemma~\ref{lem:threeisog}.

\emph{Proof of~(ii).} Let $M=\Q(\zeta_3,\sqrt[3]2)$, a Galois extension of $\Q$ with group $S_3$ generated by $r\colon\sqrt[3]2\mapsto\rho_3\sqrt[3]2$, $\rho_3\mapsto\rho_3$, and $c\colon\rho_3\mapsto\rho_3^2$, $\sqrt[3]2\mapsto\sqrt[3]2$. By~(i) and the isomorphism $E_F\cong E_{16F}$ over $M(x,y)$, the space $\mathcal N=E_F(M(x,y))\otimes\Q$ has dimension $6$ and is spanned by $P_i,[\rho_3]P_i$. For $s\in\Gal(M/\Q)$ we have $s([\rho_3]P)=[s(\rho_3)]\,s(P)$, and $r(P_i)=[\rho_3^2]P_i$ because $r(\sqrt[3]4)=\rho_3^2\sqrt[3]4$, while $c(P_i)=P_i$. Thus $\mathcal N=\bigoplus_i\Q(\rho_3)P_i$, where $r$ acts on each summand as multiplication by $\rho_3^2$ and $c$ as complex conjugation of the coefficient. The invariants of $\langle r\rangle$ are $0$, those of $\langle c\rangle$ are $\bigoplus_i\Q P_i$, and those of $S_3$ are $0$. Since $E_F(k'(x,y))\otimes\Q=\mathcal N^{\Gal(M/k')}$ for every subfield $k'\subseteq M$, the ranks are as stated. By~\eqref{eq:mordellrank} with $D=1$ and $k=\Q(\zeta_3)$, the group $E_F(\Q(\zeta_3)(x,y))$ is torsion free; here $\cX_6(k)\neq\emptyset$ because $(0,0,1)$ is a non-singular rational point of $w^6=F(x,y)$ (see~\S\ref{subsec:planes-setting}). Hence $E_F(\Q(\zeta_3)(x,y))=0$, and a decomposition $F=G^2+H^3$ would give the non-zero point $(-H,G)$.

\emph{Proof of~(iii).} Let $k=\Q(\zeta_3)$ and $B=E_{16}$ with $\beta$ as in Section~\ref{sec:number}. By~\eqref{eq:mordellrank} with $D=16$ and part~(i), $\Hom_k(A,B)^{(\sigma,\beta)}$ has rank $6$. For a sextic whose only singularities are nodes and cusps, the first Betti number of the associated six-fold cyclic multiple plane is twice the multiplicity of $t^2-t+1$ in the Alexander polynomial (\cite{Li82}, \cite{CL}; see the summary in \cite[\S1]{Kl2}); here this gives $b_1(\cX_6)=6$ and $\dim A=3$. Since $V_B=V_{B,6}$ and $e_6(B)=1$, Corollary~\ref{cor:equivhom}(i), applied over $\C$, gives $6\le\phi(6)\,e_6(A)=2e_6(A)$; as $2\dim A=6=\sum_d\phi(d)e_d(A)$, we get $e_6(A)=3$ and $V_A=V_{A,6}$. Equality holds in Corollary~\ref{cor:equivhom}(i), so for every $\lambda$ either $W_A(\lambda)=0$ or $W_B(\lambda)=H_B(\lambda)$. Since $\beta$ acts on $T_0B$ by $\zeta$, we have $W_B(\bar\zeta)=0\ne H_B(\bar\zeta)$, hence $W_A(\bar\zeta)=0$ and $W_A=H_A(\zeta)$. Thus $(A,\sigma)$ is rigid with $V_A=V_{A,6}$, which is Assumption~1.2.

Choose elements $\varphi_1,\varphi_2,\varphi_3$ of $\Hom_k(A,B)^{(\sigma,\beta)}$ forming a basis of the $\Q(\zeta)$-vector space $\Hom_k(A,B)^{(\sigma,\beta)}\otimes\Q$. The preceding dimension count and Corollary~\ref{cor:equivhom}(i) identify this space with $\Hom_{\Q(\zeta)}(V_A,V_B)$. Since $V_B$ is one-dimensional over $\Q(\zeta)$, the map $(\varphi_1,\varphi_2,\varphi_3)\colon A\to B^3$ induces an isomorphism on $V_A$, so it is an isogeny defined over $k$ and intertwining $\sigma$ with $\beta\times\beta\times\beta$. As $\End_k(B)\cong\Z[\rho_3]$ is commutative, $\End_k(A)^\sigma\otimes\Q\cong M_3(\Q(\zeta))$ has rank $18=\phi(6)\cdot3^2$, which is also the rank of $\End_{\bar k}(A)^\sigma$ by Corollary~\ref{cor:centralizer}. By Theorem~\ref{thm:planes}(b), the rank of $\Alb(\tilde\cX_6)(k(\cV_m))$ is $18m$ for every number field $k\supseteq\Q(\zeta_3)$, and by Proposition~\ref{prop:torsion} (with $e_2=e_3=0$) the group is torsion free. In this case $A^\sigma(k)=0$, so Theorem~\ref{thm:planes}(a) embeds the Mordell--Weil group into the free abelian group $(\End_k(A)^\sigma)^m$; it is therefore free.

\emph{Proof of~(iv).} Apply N\'eron's specialization theorem \cite[\S11.1]{SeMW} to $E_{16F}$ over $U$ and use~(i). The relation at $(2,3)$, where $16F(2,3)=4672$ and $Q_1=(24,136)$, $Q_2=(12,80)$, $Q_3=(8,-72)$, is checked directly.
\end{proof}

\begin{remark}
The field $\Q(\sqrt[3]2)$ enters through the coefficient $-4$ in~\eqref{eq:decomps}. By Theorem~\ref{thm:example}(ii), the equivariant homomorphisms from $A$ to the curve $y^2=x^3+1$ are defined only over $\Q(\zeta_3,\sqrt[3]{2})$, while by (iii) those from $A$ to $E_{16}$ are defined over $\Q(\zeta_3)$: in the notation of Theorem~\ref{thm:galois}, with $k=\Q(\zeta_3)$, one has $k_H=\Q(\zeta_3,\sqrt[3]{2})$ for $(B,\beta)=(E_D,\beta)$ with $D=1$, whereas $k_H=\Q(\zeta_3)$ for $(B,\beta)=(E_{16},\beta)$ and for $(B,\beta)=(A,\sigma)$.
\end{remark}

\section{Corrections to \cite{S19}, \cite{S21} and \cite{S20}}\label{sec:corrections}

We list the statements of \cite{S19}, of its corrigendum \cite{S21}, and of \cite{S20} that have to be corrected, together with the results of this paper that replace them.

\subsection*{The paper \cite{S19} and its corrigendum \cite{S21}}
\begin{enumerate}[(a)]
\item \cite[Theorem~1.1]{S19}, its amended form in \cite[(1),(2)]{S21}, and \cite[Theorem~3.3]{S20} are replaced by Theorem~\ref{thm:structure} and Corollary~\ref{cor:cyclic}. Example~\ref{ex:counter}, with $k$ a number field containing $\zeta_3$, satisfies the hypotheses of all three statements and contradicts the rank formula, both as the inequality of \cite{S19} and as the equality of \cite{S21}; if moreover $\sqrt[3]2\in k$, it also contradicts the torsion term $A[n](k)$. The error in the proof is explained in Remark~\ref{rem:S19error}. The restriction to smooth projective varieties made in \cite[(1)]{S21} is part of Setting~\ref{set:cover}.
\item \cite[Theorem~1.2]{S19} and its amended form in \cite[(3)]{S21} are replaced by Corollary~\ref{cor:superelliptic}. The torsion subgroup is $J^{\tilde\iota}(k)$ and not $J[n](k)$ (Example~\ref{ex:superelliptic}), and the rank formula holds if and only if $\tilde\iota$ is central in $\End_k(J)\otimes\Q$. The deduction of $1+\tilde\iota+\dots+\tilde\iota^{n-1}=0$ in \cite[\S5]{S19} and the assertion there about independent generators are replaced by Lemma~\ref{lem:norm} and Corollary~\ref{cor:superelliptic}.
\item Misprints in \cite[\S5]{S19}, in the notation there: on p.~335, in the description of $K=L^G$, the last generator $y_1^{s-1}y_{n-1}$ should read $y_1^{s-1}y_n$; and in~(1) on the same page, the points should read $P_1=(x_1,1)$ and $P_{i+1}=(x_{i+1},z_i/f(x_1))$ for $1\le i\le n-1$.
\end{enumerate}

\subsection*{The paper \cite{S20}}
\begin{enumerate}[(1)]
\item \cite[Theorem~3.3]{S20} is treated in item~(a) above.
\item \cite[Theorem~1.1]{S20} is replaced by Theorem~\ref{thm:planes}(b),(c): $\End_k(\Alb(\cX_n))$ must be replaced by the centralizer of $\bar\tau_*$, and the torsion term $\Alb(\cX_n)[m](k)$ by $A^\sigma(k)$.
\item The deduction of $\id+\tilde\gamma+\dots+\tilde\gamma^{n-1}=0$ in \cite[\S5]{S20} is replaced by Lemma~\ref{lem:norm}.
\item \cite[Theorem~1.3]{S20} is replaced by Corollary~\ref{cor:number}. Under Assumption~1.2 the rank is at most $m\cdot4d_n^2/\phi(n)$ for every $n\ge3$; the value $2d_n^2$ given for $n=5,10$ is incorrect when $d_n>0$. The equalities hold only if the relevant endomorphisms are defined over $k$ (the decompositions in \cite[Prop.~4.2]{S20} are over $\C$). The restriction on $n$ is unnecessary.
\item The argument for $n=8$ in \cite[\S6]{S20} is incorrect. The curve $y^2=x^3-x^2-3x+1$, denoted $E_2$ in \cite{S20}, has $c_4=160$ and discriminant $2^6\cdot37$, so $j(E_2)=64000/37\notin\Z$ and $E_2$ has no complex multiplication; in particular $J(C_2)$ is not isogenous to $E_1\times E_2$. By Remark~\ref{rem:n812}(2), $J(C_2)$ is isogenous over $\C$ to the square of an elliptic curve with complex multiplication by $\Z[\sqrt{-2}]$; for instance the curve $y^2=x^3+x^2-3x+1$, denoted $E_1$ in \cite{S20}, which has $j$-invariant $8000$. Also, $M_{2n_2}(\Z[i])$ has rank $8n_2^2$, not $4n_2^2$. Under Assumption~1.2, one of $n_1,n_2$ in \cite[Prop.~4.2(v),(vi)]{S20} is zero (Remark~\ref{rem:n812}).
\item \cite[Theorem~1.4]{S20} is replaced by the discussion after Corollary~\ref{cor:number}. For the dual of the Fermat cubic the conclusion $2md_6^2=18m$ holds over every number field containing $\Q(\zeta_3)$ (Theorem~\ref{thm:example}(iii)). For the curve in \cite[Theorem~1.4(ii)]{S20}, Assumption~1.2 was not verified.
\item \cite[Remark~1.5(1)]{S20} holds in the form of Remark~\ref{rem:higherdim}. \cite[Remark~1.5(2)]{S20} is justified only for $m=1$ (Proposition~\ref{prop:specialization} and Remark~\ref{rem:mge2}).
\item Misprints in \cite{S20}:
\begin{itemize}
\item p.~861: ``$r$ is a positive number'' should read ``$r$ is a non-negative integer''.
\item p.~863: ``the function fields of $U_m$ and $V_m$ respectively'' should read ``the function fields of $V_m$ and $U_m$, respectively''; ``the function field of $W_m$'' should read ``the function field of $V_m$''.
\item p.~864, Theorem~1.4(ii), and p.~875: the affine equation should be $f(x,y)=F(1,x,y)=x^3-3xy(y^3-8)+2(y^6+20y^3-8)$, with $F$ as in \cite[Prop.~4.4(ii)]{S20}; the sign of $20y^3$ is misprinted.
\item p.~869: ``$e\ge n/r$'' should read ``$e\ge r/n$'', and the weighted projective space $\PP^3_{(1,1,1,n_0)}$ should be $\PP^3_{(1,1,1,e)}$, so that both sides of $u_3^n=u_0^{n_0}F(u_0,u_1,u_2)$ have degree $ne=n_0+r$.
\item p.~871, Prop.~4.4(ii): the factor $(u_2^3-8)$ should read $(u_2^3-8u_0^3)$.
\item p.~874: $\End_k(\Alb(\cX_5))$ should read $\End_k(\Alb(\cX_8))$ (three times).
\item p.~875, proof of Theorem~1.4: ``$2md_6$'' should read ``$2md_6^2$''.
\item p.~875, reference~[9]: the year of the English translation is $1995$, not $1955$.
\end{itemize}
\end{enumerate}


\begin{thebibliography}{99}

\bibitem{BS} \textsc{A.~Borel and J.-P.~Serre}, \textit{Th\'eor\`emes de finitude en cohomologie galoisienne}, Comment. Math. Helv. \textbf{39} (1964), 111--164.

\bibitem{BLR} \textsc{S.~Bosch, W.~L\"utkebohmert and M.~Raynaud}, \textit{N\'eron models}, Ergebnisse der Mathematik und ihrer Grenzgebiete (3) \textbf{21}, Springer-Verlag, Berlin, 1990.

\bibitem{CL} \textsc{J.~I.~Cogolludo-Agust\'in and A.~Libgober}, \textit{Mordell--Weil groups of elliptic threefolds and the Alexander module of plane curves}, J. Reine Angew. Math. \textbf{697} (2014), 15--55.

\bibitem{De} \textsc{A.~Degtyarev}, \textit{Irreducible plane sextics with large fundamental groups}, J. Math. Soc. Japan \textbf{61} (2009), no.~4, 1131--1169.

\bibitem{Es} \textsc{H.~Esnault}, \textit{Fibre de Milnor d'un c\^one sur une courbe plane singuli\`ere}, Invent. Math. \textbf{68} (1982), no.~3, 477--496.

\bibitem{Ha} \textsc{F.~Hazama}, \textit{Rational points on certain abelian varieties over function fields}, J. Number Theory \textbf{50} (1995), no.~2, 278--285.

\bibitem{HK} \textsc{K.~Hulek and R.~Kloosterman}, \textit{Calculating the Mordell--Weil rank of elliptic threefolds and the cohomology of singular hypersurfaces}, Ann. Inst. Fourier (Grenoble) \textbf{61} (2011), no.~3, 1133--1179.

\bibitem{Kl1} \textsc{R.~Kloosterman}, \textit{Cuspidal plane curves, syzygies and a bound on the MW-rank}, J. Algebra \textbf{375} (2013), 216--234.

\bibitem{Kl2} \textsc{R.~Kloosterman}, \textit{Mordell--Weil lattices and toric decompositions of plane curves}, Math. Ann. \textbf{367} (2017), no.~1--2, 755--783.

\bibitem{Ku} \textsc{Vik.~S.~Kulikov}, \textit{On plane algebraic curves of positive Albanese dimension}, Izv. Ross. Akad. Nauk Ser. Mat. \textbf{59} (1995), no.~6, 75--94; English transl., Izv. Math. \textbf{59} (1995), no.~6, 1173--1192.

\bibitem{Li82} \textsc{A.~Libgober}, \textit{Alexander polynomial of plane algebraic curves and cyclic multiple planes}, Duke Math. J. \textbf{49} (1982), no.~4, 833--851.

\bibitem{Li12} \textsc{A.~Libgober}, \textit{Factors of Jacobians and isotrivial elliptic surfaces}, J. Singul. \textbf{5} (2012), 115--123.

\bibitem{Li13} \textsc{A.~Libgober}, \textit{On Mordell--Weil groups of isotrivial abelian varieties over function fields}, Math. Ann. \textbf{357} (2013), no.~2, 605--629.

\bibitem{Mi} \textsc{J.~S.~Milne}, \textit{Abelian varieties}, in: Arithmetic Geometry (Storrs, Conn., 1984), G.~Cornell and J.~H.~Silverman (eds.), Springer, New York, 1986, 103--150.

\bibitem{Ne} \textsc{A.~N\'eron}, \textit{Probl\`emes arithm\'etiques et g\'eom\'etriques rattach\'es \`a la notion de rang d'une courbe alg\'ebrique dans un corps}, Bull. Soc. Math. France \textbf{80} (1952), 101--166.

\bibitem{Ni} \textsc{H.~Nishimura}, \textit{Some remark on rational points}, Mem. Coll. Sci. Univ. Kyoto Ser.~A Math. \textbf{29} (1955), 189--192.

\bibitem{S19} \textsc{S.~Salami}, \textit{The rational points on certain Abelian varieties over function fields}, J. Number Theory \textbf{195} (2019), 330--337.

\bibitem{S20} \textsc{S.~Salami}, \textit{Twists of the Albanese varieties of cyclic multiple planes with large ranks over higher dimension function fields}, J. Th\'eor. Nombres Bordeaux \textbf{32} (2020), no.~3, 861--876.

\bibitem{S21} \textsc{S.~Salami}, \textit{Corrigendum to ``The rational points on certain Abelian varieties over function fields'' [J. Number Theory 195 (2019) 330--337]}, J. Number Theory \textbf{227} (2021), 306--307.

\bibitem{SeAlb} \textsc{J.-P.~Serre}, \textit{Morphismes universels et vari\'et\'es d'Albanese}, S\'eminaire Claude Chevalley, tome~\textbf{4} (1958/59), expos\'e~10.

\bibitem{SeTopics} \textsc{J.-P.~Serre}, \textit{Topics in Galois theory}, Research Notes in Mathematics~\textbf{1}, Jones and Bartlett, Boston, 1992.

\bibitem{SeMW} \textsc{J.-P.~Serre}, \textit{Lectures on the Mordell--Weil theorem}, 3rd ed., Aspects of Mathematics~\textbf{E15}, Vieweg, Braunschweig, 1997.

\bibitem{Sh} \textsc{T.~Shioda}, \textit{On the Mordell--Weil lattices}, Comment. Math. Univ. St. Paul. \textbf{39} (1990), 211--240.

\bibitem{Si} \textsc{J.~H.~Silverman}, \textit{Heights and the specialization map for families of abelian varieties}, J. Reine Angew. Math. \textbf{342} (1983), 197--211.

\bibitem{Ta} \textsc{J.~Tate}, \textit{Algorithm for determining the type of a singular fiber in an elliptic pencil}, in: Modular Functions of One Variable IV, Lecture Notes in Math. \textbf{476}, Springer, Berlin, 1975, 33--52.

\bibitem{TS} \textsc{J.~T.~Tate and I.~R.~Shafarevich}, \textit{The rank of elliptic curves}, Dokl. Akad. Nauk SSSR \textbf{175} (1967), 770--773.

\bibitem{Ul} \textsc{D.~Ulmer}, \textit{On Mordell--Weil groups of Jacobians over function fields}, J. Inst. Math. Jussieu \textbf{12} (2013), no.~1, 1--29.

\bibitem{Wa} \textsc{W.~B.~Wang}, \textit{On the twist of abelian varieties defined by the Galois extension of prime degree}, J. Algebra \textbf{163} (1994), no.~3, 813--818.

\bibitem{Ya} \textsc{H.~Yamagishi}, \textit{A unified method of construction of elliptic curves with high Mordell--Weil rank}, Pacific J. Math. \textbf{191} (1999), no.~1, 507--524.

\end{thebibliography}
\end{document}